\documentclass{article}

\usepackage{PRIMEarxiv}

\usepackage[utf8]{inputenc} 
\usepackage[T1]{fontenc}    
\usepackage{hyperref}       
\usepackage{url}            
\usepackage{booktabs}       
\usepackage{amsfonts}       
\usepackage{nicefrac}       
\usepackage{microtype}      
\usepackage{lipsum}
\usepackage{fancyhdr}       
\usepackage{graphicx}       
\graphicspath{{media/}}     

\usepackage{natbib}
\bibpunct[, ]{(}{)}{,}{a}{}{,}%
\usepackage{subfigure}
\usepackage{lmodern} 
\usepackage{array}
\usepackage{booktabs}
\usepackage{array, threeparttable}
\usepackage{tikz}
\usepackage{pifont}
\usepackage{graphicx} 
\usepackage{float} 
\usepackage{makecell}
\usepackage{amsmath}
\usepackage{algorithm}
\usepackage{algorithmic}
\usepackage{graphicx} 
\usepackage{multirow}
\usepackage{subfigure}
\usepackage{color}
\usepackage{threeparttable}
\usepackage{tabularx}
\usepackage[figuresright]{rotating}
\usepackage{caption}
\usepackage{adjustbox}  
\usepackage{enumerate}
\usepackage{setspace}
\usepackage[resetlabels]{multibib}
\usepackage{hyperref}
\newcites{sec}{Reference}
\newcommand{\deq}{\overset{\mathrm{def}}{=}}
\title{City-LEO: Toward Transparent City Management Using LLM with End-to-End Optimization}

\author{
  Zihao Jiao\\
  School of Computer Science and Artificial Intelligence \\
  Beijing Technology and Business University \\
  Beijing, China\\
  \And
  Mengyi Sha\\
  Department of Industrial Engineering\\ 
  Tsinghua University\\
  Beijing, China\\
    \And
  Haoyu Zhang\\
  School of Computer Science and Artificial Intelligence \\
  Beijing Technology and Business University \\
  Beijing, China\\
    \And
  Xinyu Jiang\\
  Department of Industrial Engineering\\ 
  Tsinghua University\\
  Beijing, China\\
    \And
  Wei Qi\\
  Department of Industrial Engineering\\ 
  Tsinghua University\\
  Beijing, China\\
  \texttt{qiwei.0216@gmail.com}
}

\begin{document}
\maketitle

\begin{abstract}
Existing operations research (OR) models and tools play indispensable roles in smart-city operations. However, the adoption of conventional OR tools confronts some challenges: On one hand, the highly intricate nature of urban operations introduces difficulties in accurately modeling the problems in city management. On the other hand, the inherent intransparency of conventional OR models, coupled with the requirement for specialized OR expertise, also sets barriers in the development and implementation processes. Amidst these challenges, Large Language Model (LLM) is rapidly evolving and offers a promising path to bridge the gap between OR theory and practical implementation. Nevertheless, existing LLM-based OR techniques exhibit several limitations in terms of relevance and accuracy. To alleviate the limitations, in this paper we make an early attempt to propose an innovative LLM-embedded agent, City-LEO. The agent aims to conduct human-like decision process and deliver more relevant and accurate solutions to tackle city management problems. Specifically, to enhance computational tractability and relevance to users’ queries, based on a pre-specified full-scale optimization model, City-LEO leverages LLM’s logical reasoning capability on prior knowledge to effectively scope down the users’ requirements. In addition, to promote transparency and account for environmental features, City-LEO incorporates an End-to-End (E2E) framework to synergize prediction and optimization process. The feature-to-decision framework is conducive to coping with environmental uncertainties and involving social concerns. In the case study, we employ City-LEO in the operations management of an electric-bike sharing (EBS) system. We design relevance tests and accuracy tests to evaluate the performance of City-LEO, and the numerical results demonstrate that City-LEO has superior performance when benchmarked against the full-scale optimization problem. With less computational time, City-LEO generates more relevant solutions to the users' requirements without significantly compromising accuracy. In a broader sense, City-LEO showcases the potential of developing LLM-embedded OR tools for smart-city operations management.
\end{abstract}


\section{Introduction}\label{sec:Intro}
Operations research (OR) models and tools, such as commercial solvers including Gurobi and Cplex, are being used extensively in smart-city transportation management to prescribe planning and operational decisions \citep{qi2019smart, zhang2022dynamic, gayialis2022city, chen2023RSOME}. For urban bus transit, bike-sharing and ride-hailing systems, OR models and tools facilitate the decision-making processes and generate planning and scheduling strategies to promote operational efficiency. 
Nevertheless, as \citep{bettencourt2024recent} in \textit{Science} points out, the urban systems encompass a multitude of interrelated or uncertain factors, and are usually highly complex to model accurately. The highly intricate nature of urban operations, coupled with the heightened need for specialized OR applications, results in a significant disparity between theory and practical implementation. 
   The adoption of OR tools also confronts challenges that stem from inherent \emph{intransparency} in development and implementation processes. Generally, city operators possess a deep understanding of urban operations mechanisms but lack a fundamental understanding of OR expertise. Conversely, the developers of OR tools have a solid grasp of OR methodologies, but often lack knowledge of the practical operational principles of urban systems. Moreover, conventional OR tools cause information asymmetry between the city operators and the residents. The city operators hold a dominant position in the decision-making process and control the OR tools, while the process appears to be a black box for the residents. 

Given these challenges at hand, large language models (LLMs) hold promise to empower city operators to wield OR models and solvers toward enabling more transparent city management. 
The boom of LLMs, such as GPT and Ollama, marks a breakthrough in artificial intelligence.  Based on deep learning techniques, LLMs offer a path to meaningful utilization of the growing amounts of municipal data and understanding, summarizing, generating, and predicting new content. LLMs have recently penetrated into an array of application fields. Several vertical fields of LLMs, such as the legal LLM ChatLaw and the financial LLM XuanYuan 2.0 \citep{zhang2023xuanyuan}, have emerged to prescribe intelligent decisions in their respective fields. As a promising technique for assisting humans in multiple tasks, LLMs have also been tried in OR applications. In practice, there have been two main pathways for implementing LLMs:
\begin{itemize}
    \item \emph{LLMs as optimizers}. This technique leverages the in-context learning LLM \citep{dong2022survey} to directly address an optimization problem conversationally without invoking an optimization model. By integrating chain of thought (CoT) prompts \citep{kojima2022large}, the LLM's reasoning abilities can be harnessed to directly tackle and resolve an optimization problem. LLMs as optimizers do not require extensive knowledge of specific details or specialized OR expertise. The technique only needs a few prompts (e.g., \emph{Let's think step by step} \citep{wei2022chain}) to provide solutions for user-defined problems.
    \item \emph{LLM agents}. In order to enhance the transparency and accessibility of the decision-making process, another straightforward strategy would be to develop an LLM agent that is capable of generating optimization models (in the form of callable APIs). Some recent preliminary works propose to integrate OR optimization techniques with LLM into an \textit{Agent}, and these LLM-based agents mainly focus on deterministic mixed integer programming (MIP) models \citep{kasneci2023chatgpt, li2023large,ahmaditeshnizi2023optimus}.
\end{itemize}

\begin{figure}[] 
\centering 
\includegraphics[width=0.65\textwidth]{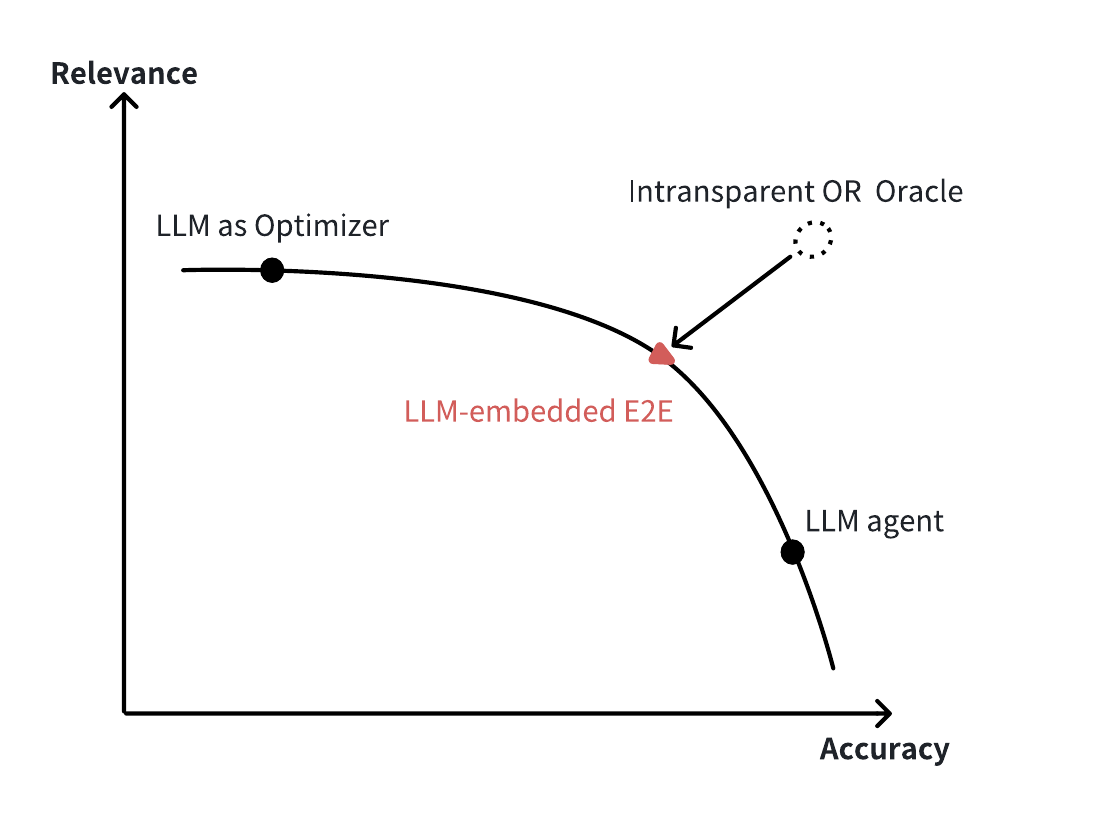} 
\caption{Dilemma between Relevance and Accuracy of LLM-based OR Tools in City Management} 
\label{Fig.dilemma} 
\end{figure}

Nevertheless, existing LLM-based OR techniques exhibit limitations in terms of \emph{relevance} and \emph{accuracy}. The extent to which we can extrapolate from these LLM capabilities to the intricate realm of city management remains unclear. Although \emph{LLMs as optimizers} negate the need for OR expertise, guaranteeing the \emph{accuracy} of solutions without the aid of rigorous modeling and solvers is challenging. For example, the experiments of LLM optimizers show that there is still an 11\% gap between the achieved solution and the optimal solution in a 50-node TSP problem \citep{yang2023large}. On the other hand, \emph{LLM agents} are generally restricted to dealing with simple and well-defined problems. However, city operators usually confront highly intricate issues, and may struggle to accurately match suitable models and code with vague and macro-level queries, which may result in \emph{irrelevant} optimization model recommendations. Consequently, agents associated with callable MIPs (e.g., \citep{ahmaditeshnizi2023optimus} and \citep{li2023large}) may lack the flexibility to generate adequate functionalities and accommodate city operators' diverse requests. Moreover, agents with callable MIP often overlook environmental uncertainties, underutilize historical data and involve limited parameters and decisions, which hinders their applicability and effectiveness in addressing complex urban operational problems. On the whole, as Figure \ref{Fig.dilemma} illustrates, there is a dilemma between the relevance and accuracy of LLM-based OR tools in city management. The dilemma lies in the difficulty of simultaneously achieving complex problem-solving and optimality. Specifically, in most cases, it is difficult to prescribe a model in OR applications that can match users' diverse queries based on large-scale complex parameters and decision types and then use a solver to get efficient and exact solutions. 
The limitations of \emph{LLMs as optimizers} and \emph{LLMs agents} urge us to identify a sweet middle ground to get close to the intransparent (and often intricate) OR Oracle (see, Figure \ref{Fig.dilemma}) without much compromising the relevance or accuracy. Here, we are motivated to propose City-LEO, an {\bf L}LM-based agent joint with {\bf E}2E {\bf O}ptimization, capable of providing reasonably relevant and accurate operational strategies to city operators through transparent interactions: 
\begin{itemize}
\item To ensure accuracy and relevance of the LLM-embedded optimization tool, we develop a full-scale optimization model first. The pre-defined model involves more constraints and covariates (e.g. features and decision variables), and allows for a nuanced characterization of the complex operations of a specific urban system. In this way, the comprehensive model incorporates domain knowledge and serves as the bedrock of our proposed framework.

\item To further enhance relevance and computational tractability, we propose a human-like decision-making process. The process allows City-LEO to leverage LLM's reasoning capability based on prior knowledge, so as to effectively \emph{scope down} the pre-specified full-scale optimization model based on users' requirements. Such problem-scoping operations allow for focusing on a careful choice of appropriate parameters and decisions from a diverse set of potential candidates for problem modeling to minimize irrelevance and computational expenditure. Additionally, the daily operational goals proposed by city operators usually focus on \emph{neighborhood-level} problems, which typically do not necessitate redundant full-scale optimization models as seen in traditional OR applications. Thus, City-LEO is conducive to effectively scoping down the size of the large-scale optimization problem to accommodate users' diverse requirements, while lowering the expertise barrier for city operators in urban governance.

\item In order to further enhance transparency and account for environmental uncertainties and contextual features, we extend from the literature an E2E framework that outputs decisions directly from input covariates. Specially, the E2E framework integrates prediction with optimization processes, and conducts MIP optimization over a Random Forest (RF) objective function with general polyhedral constraints. On the one hand, the tool entails adopting a large number of features into the optimization model, which has the potential to increase the degree of relevance to the user's requirements and then mitigate decision-making bias. On the other hand, random forest is powerful to depict the complex and unknown nonlinear relationships between features and macro-level indicators, and the hierarchical tree structure can reflect the decision paths in random forest and makes the E2E framework more interpretable and transparent.
\end{itemize}
Taken together, we summarize the contributions of our work as follows: 
\begin{itemize}
    \item[(1)] We make an early attempt to propose an innovative LLM-embedded agent, City-LEO, with the aim of delivering more \emph{relevant} and \emph{accurate} solutions to tackle city management problems. The agent also enhances the accessibility of OR tools through conversational interactions, and promotes the transparency and efficiency of the decision-making processes.
    \item[(2)] We expand the nascent agent framework of synergizing LLMs and OR in a methodology-wise approach. In particular, (i) City-LEO is conducive to effectively scoping down the pre-defined large-scale optimization problem to accommodate users' diverse requirements. The agent leverages the LLM's logical reasoning ability based on prior knowledge to redefine the problem and conduct human-like decision-making processes. The scope-down policy not only enhances the computational tractability of the optimization model, but also elevates the relevance of the solutions in addressing city management problems. (ii) We integrate an E2E modeling framework with a pre-trained random forest to cope with city management problems in an uncertain operational environment. Different from existing MIP agents, the E2E framework integrates estimation with optimization, and offers an opportunity to generate decisions directly from feature data, which is more familiar to residents. The feature-to-decision framework promotes data utilization \citep{qi2023e2e}, and also leads to more transparent and responsible decision processes by incorporating indicators to reflect social concerns. 
    \item[(3)] We numerically test City-LEO in the operations management of an urban electric bike sharing (EBS) system. Specifically, the agent focuses on optimizing the e-bike rebalancing and battery swapping operations to mitigate the imbalance between supply and customers' demands. Our case study demonstrates that the City-LEO has superior performance in terms of accuracy and relevance. By benchmarking against the full-scale MIP model, the numerical results validate the effectiveness of City-LEO in terms of computational time and the quality of the solutions. Moreover, the results demonstrate that our proposed approach holds promise of using LLM to amplify the value of OR models in addressing complex smart-city operations problems. 
    \end{itemize}

The remainder of this paper is organized as follows. Section \ref{sec:LR} reviews the related literature. Section \ref{sec:agent} provides a formal introduction to the structure of the LLM agent. Section \ref{sec:e2e} presents the E2E models for the EBS system management. Section \ref{sec:CS} discusses numerical results. Section \ref{sec:con} concludes the paper. Notation tables, methodology review, additional experimental settings and sample of prompts are available in the online appendices.
We also publish our data and code source related to the E2E model, City-LEO agent, and necessary materials files to  Hugging-face (please see \href{https://huggingface.co/datasets/BTBUzhy/JOC}{City-LEO})
for the convenience of readers.

\section{Literature Review}\label{sec:LR}

Synergizing OR and LLM on tackling complex urban management problems is an important and intriguing paradigm for future smart cities, but so far has not been examined in the literature to our knowledge. Existing literature has studied both utilizing the reasoning ability of LLM to solve optimization problems heuristically and integrating LLM with professional tools by constructing sophisticated structures for agents,  yet with little attention on LLM-embedded optimization. We review the related literature in multiple academic fields, and then highlight the distinction of our proposed methodology.
To familiarize the OR audience with a more general background on using LLMs for optimization, we provide an extended literature survey on the recent development of In-context learning in Section \ref{app:ICL}. 

\subsection{LLM Agents}
LLM agents \citep{wu2023autogen,zhao2024expel}, which integrate externally callable domain knowledge and flexible functions, are capable of mitigating the low accuracy and token-limitation (as demonstrated in Section \ref{app:ICL}) issues associated with in-context learning. 
Typically, LLM agents utilize external functions to generate various extensions, which can be classified as \emph{LLM agents with function calling} and \emph{Autonomous language agents} (ALA, see Section \ref{app:ALA} for details). 

\textbf{LLM Agents with Function Calling.} The first approach entails directly utilizing pre-packaged external functions, referred to as "function calling" (LLM-FC). To enhance the accuracy of the responses, LLM-FC leverages the reasoning capabilities of LLM to match the user's query with domain knowledge (e.g., open data source, commercial solvers like Gurobi and Cplex), which is pre-defined as structural functions in an external API \citep{ge2024openagi}. In the context of optimization, the main concept underlying LLM-FC is not to substitute optimization technology with LLMs but rather to employ optimization solvers in tandem with LLMs. The agent utilizes a specific and specialized template to identify users' queries and generate exact solutions to offer responses. For instance, through commercial solver agents \citep{ahmaditeshnizi2023optimus}, users can adopt natural language to invoke code generation commands in the solver. \citep{li2023large} proposes agents and a response architecture to cope with "what if" questions. The agent modifies the input to the optimization solver in a suitable manner and subsequently re-executes the solver to generate a solution.

Acting as a "translator" between users and expertise OR tools, the  paradigm associated with LLM agents confronts with professional barriers yet:  Typically, city administrators who lack expertise in optimization theory struggle to accurately match suitable models and precise code with vague and macro-level queries, resulting in ineffective optimization recommendations.  
For example, prior to proceeding, LLM-FC  requires the user to provide an exceedingly accurate and detailed depiction of the issue, including quantifiable objectives, key parameters (e.g. unit cost and resources) and constraints to represent the problem setting \citep{ahmaditeshnizi2023optimus,li2023synthesizing}. 

However, in practice, it is challenging for users to provide  precise and detailed description of the problem. Thus, the issue sets barriers to extracting useful information from users' queries and constructing "accurate" and "feasible" models for LLM-FC, which may result in inefficient responses and weakening the model's potential for practical use. 
For instance, \citep{ahmaditeshnizi2023optimus} shows that the success rate (i.e., the ratio of outputs that satisfies constraints and can find the optimal solution) is less than 80\% based on GPT-4.0 function calling.
Hence, the subsequent limitations of employing such agents cause possible discrepancies between users' desired responses and the generated code. Consequently, MIP-embedded agents possibly prescribe wrong and opaque (if not infeasible) decisions for citywide management. 
We summarize the limitations of the existing approaches as follows. Area-specific MIP agents (e.g. LLM-FC and ALA), relying on deterministic assumptions, are unable to provide a diverse range of solutions that effectively meet the users' requirements. The efficacy of the agents are limited by the pre-defined parameters and decision variables. Hence, it is imperative to develop decision-support LLM agents that can effectively address users' queries and comprehensively consider the attributes of urban systems and environments.
\subsection{End-to-End Optimization}
LLM agents have some limitations in dealing with various practical problems in uncertain environments. (i) If a user's query does not cover some key parameters or conditions of the MIP model, the code generation is either incomplete or ambiguous. (ii) Users' open-ended queries consistently cause variable problem settings, which poses challenges to deterministic models. Moreover, deterministic MIP optimization models incur biased solutions due to inefficient estimation of uncertainties. In contrast, E2E optimization offers a promising way for LLM-assisted optimization to generate relevant and accurate solutions.
E2E optimization refers to outputting decisions directly from the input covariates. The feature-to-decision mapping integrates prediction and optimization process.

We identify two branches of research on feature-to-decision mapping: (i) One-step mapping is an integrated version of learning and optimization to search for the best-performing decisions. 
\citep{qi2023e2e} involves deep learning models in the E2E framework to generate the multi-period inventory replenishment strategies directly from the input features. 
(ii) Two-step mapping means sequential learning and optimization process, which initially utilizes given covariates and trained model to estimate the conditional distribution of the uncertain parameters, and then tackles an associated optimization problem to derive the optimal solution \citep{elmachtoub2022smart}. For instance, \citep{Biggs_2022} propose to reformulate random forests along with operational polyhedral constraints into 0-1 mixed integer programs.
The prediction in this approach ignores the properties of subsequent optimization models, which may result in error propagation between estimation and optimization and suboptimal solutions \citep{qi2022integrating}.

Therefore, in our proposed LLM-based optimization tool, we adopt E2E optimization to address open-ended queries about EBS management for the following reasons: First, extending the approach of \citep{Biggs_2022}, utilizing pre-trained random forests can avoid repetitive training in dealing with diverse user queries. Second, the objectives can include indicators (e.g., traffic congestion metrics) related to social concerns and contribute to mitigating negative social externalities. For city administrators, the tool based on E2E optimization is user-friendly and does not require extensive specialized knowledge.

\subsection{In-context Learning}\label{app:ICL}

One topic that aligns with our research is In-context learning. Research in this area is mostly aimed at promoting LLMs' abilities to solve area-specific problems (such as mixed integer programming in our context) through \emph{Prompt Engineering}, represented as in-context learning methods \citep{dong2022survey}, which include the zero-shot chain of thought (CoT) \citep{kojima2022large}, few-shot CoT \citep{wei2022chain}, and Auto-CoT \citep{zhang2022automatic}.
Specifically, zero-shot CoT refers to the capacity of LLMs to perform logical reasoning without any prior training examples, whereas few-shot CoT entails offering the model with a limited number of examples to direct its reasoning.
The auto-CoT samples questions and develops chains of reasoning to construct demonstrations. One of the relevant research \citep{yang2023large} leverages the capabilities of LLMs as optimizers, designs a meta-prompt, and solves optimization problems interpreted from natural language. Nevertheless, the results show a poor level of accuracy in solving the optimization problem. As an illustration, the LLM optimizer experiment still shows an 11\% optimality gap from the best solution in a 50-node TSP problem, as determined using oracle solutions. 

The low accuracy observed in the aforementioned CoT methods is caused by the model's unstable reasoning results, resulting in hallucinatory phenomena such as factual inaccuracies and logical errors.  \citep{chu2023survey} attributes the issue to the absence of high-quality artificial annotations and limited reasoning paths. To mitigate hallucinatories, \citep{yao2023tree} proposed a novel framework called the tree of thought (ToT), which enables the LLM to take into account various reasoning paths. ToT greatly enhances the interpretability and reasoning capabilities of the reasoning process by adding more prompts. However, limited by capacity and resources, LLMs do not have enough tokens to deal with longer prompts. TOT is subject to significantly more constraints when choosing activities and generating tailored text recommendations for each query.


\subsection{Autonomous Language Agent}\label{app:ALA}
ALA has the ability to automatically handle objective-oriented multi-step tasks \citep{park2023generative}, and enables autonomous problem-solving, allowing it to handle specific objectives and tasks independently rather than simply responding to users' queries.

In this vein, \citep{yao2023retroformer} utilizes verbal "feedback", specifically self-reflection, to assist agents in acquiring knowledge from previous instances of failure. Furthermore, the feedback can be defined as gradient reduction of optimization objectives, thereby contributing to solving MIP model. For instance, \citep{yao2022react} propose a reinforced large language agent by adopting environment-specific rewards. The agent utilizes the rewards from various contexts and tasks to fine-tune a pre-trained LLM. The LLM is then used to improve the language agent prompt by summarizing the main reason for the previous unsuccessful attempts and suggest action plans. Nevertheless, the ALA continues to exhibit low accuracy and significant token expenses due to the presence of a large-scale ICL prompts. 

\section{The City-LEO Agent}\label{sec:agent}
This section illustrates the task flow framework and an EBS case for our proposed City-LEO agent. 
Through the \textcircled{1} \textit{Problem matcher}, the agent matches (i.e., first-stage scoping down) users' natural language queries $Q$  (See Appendix \ref{app:notation} for all notations) with appropriate area-specific agents (i.e., callable tool APIs). It consists of the following three critical components: \textcircled{2} The \textit{Query-relevant objective generator} is the process of generating an Query-relevant objective (QR-obj) function $f(\textbf{x};\textbf{w})\deq\textbf{LLM}(Q;P^{\text{IG}}): \mathbb{R}^{|\mathcal{X}|\times |\mathcal{W}|}\rightarrow \mathbb{R}_+$ for input $Q$ by incorporating pre-defined prompts $P^{\text{IG}}$,  where $\mathcal{X},\mathcal{W}$  denote the feasible domain of decision variables $\textbf{x}$ and set of parameters  $\textbf{w}$, respectively. Function $f(\textbf{x};\textbf{w})$ serves as a secondary objective in E2E optimization problem $\min_{\textbf{x}\in \mathcal{X}} g(\textbf{x};\textbf{w})$.
\textcircled{3} \textit{Problem tailor} is adopted for second-stage scoping down the feasible domain $\mathcal{X}$ of decision variables $\textbf{x}$, particularly searching for \emph{parameterized variables} $\textbf{x}^{\prime}\in \mathcal{X}$ in the E2E optimization programming $\min_{\textbf{x}\in \mathcal{X}} g(\textbf{x};\textbf{w})$.
Parameterized variables refer to the decision variables that are less relevant to the user's query, based on an inference of historical information. And we assign historical values to parameterized variables in the form of equation constraints.
\textcircled{4} The \textit{E2E model} comprises the pre-trained random forest model, along with the associated MIP model. 
The agent incorporates a code safeguard that performs iterative examinations based on the optimization language rule of Gurobi. Figure \ref{Fig.framework} provides a comprehensive illustration of the agent framework:

\textbf{Step 1: Problem Matcher}. The user (e.g., the city operators, the EBS system managers, etc.) submits a natural language query $Q$  based on a certain target in a specific area (e.g., \textit{How to decrease the proportion of shared e-bikes with low SOC in parking spots 5, 6, and 7?}). Such a query is subsequently processed by the \textit{problem matching agent} (see Figure \ref{Fig.framework} \textcircled{1}) to determine suitable area-specific agents (e.g., EBS, bus line analysis, or other city management area agents). Such area-specific agents encompass QR-obj generators, problem tailors, the E2E models, and response-prompters that are specifically associated with the given problem.\\
\begingroup
\setlength\abovedisplayskip{0pt}
\setlength\belowdisplayskip{0pt}
\setlength\abovedisplayshortskip{0pt}
\setlength\belowdisplayshortskip{0pt}
\begin{figure}[]
\centering 
\includegraphics[width=0.85\textwidth]{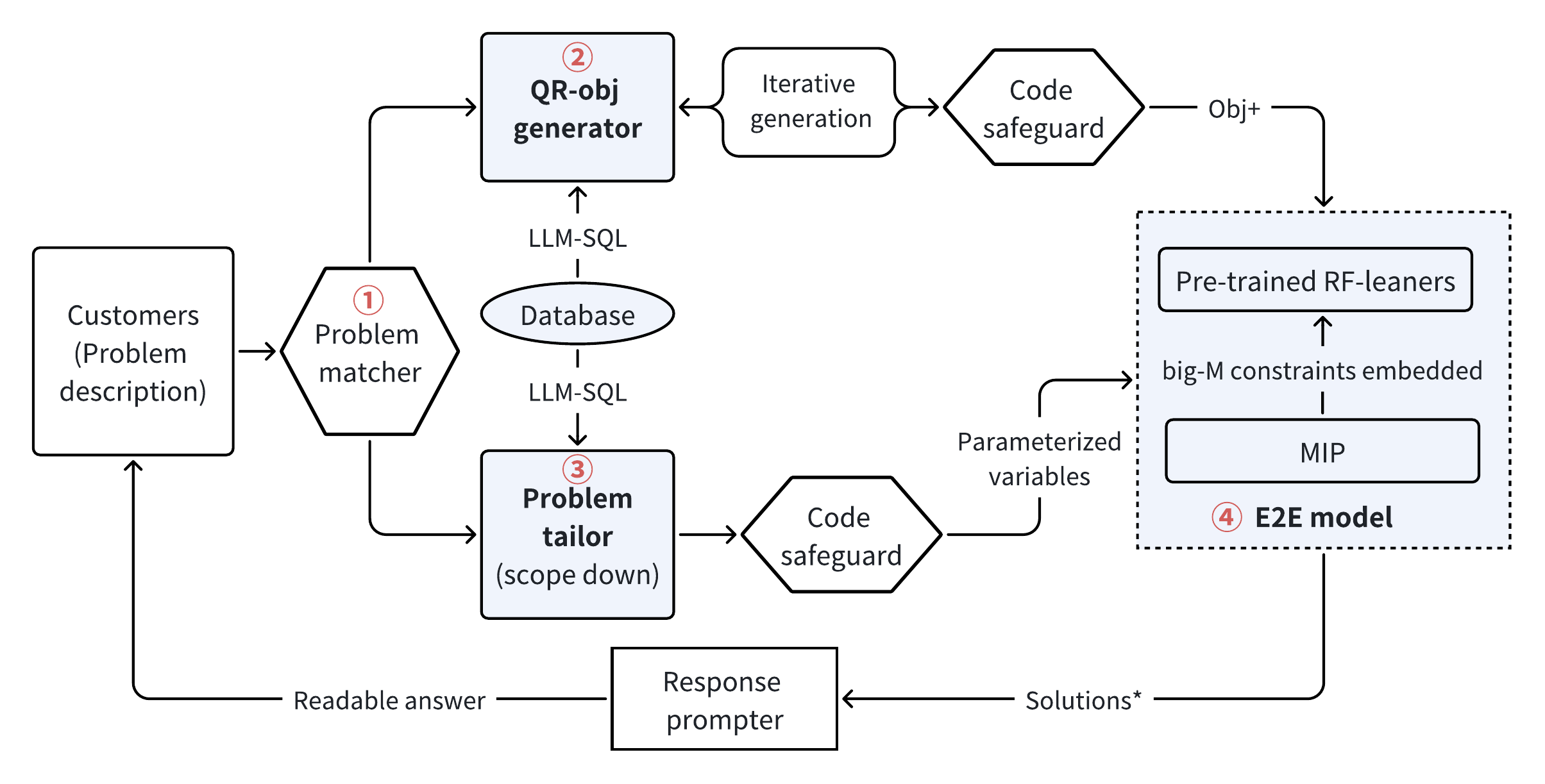} 
\caption{Framework of the City-LEO Agent} 
\label{Fig.framework} 
\end{figure}
\endgroup

\textbf{Step 2: QR-obj Generator}.  
The user query is subsequently processed by the QR-obj generator within area-specific agents that have been matched. 
The QR-obj generator (see Figure \ref{Fig.framework} \textcircled{2}) formulates natural language queries $Q$ into a convex QR-obj function $f(\textbf{x};\textbf{w})$ of decisions $\textbf{x}$ (e.g., number of dispatched e-bikes or number of swapped depleted batteries) and parameters $\textbf{w}$ in the E2E programming $\min_{\textbf{x}\in \mathcal{X}} g(\textbf{x};\textbf{w})$. 
The QR-obj $f(\textbf{x};\textbf{w})$ is then incorporated into the objective of MIP for re-optimization. 
In addition, the retrievable historical operational databases $\mathcal{D}$ (in the form of a callable API with LLM-generated structured query language (SQL)) provide adequate operational information for QR-obj generator. It further strengthen the relevancy and significance of the produced QR-obj to users' query.
Following that, a \emph{Code safeguard} will iteratively examine the output code of $f(\textbf{x};\textbf{w})$ to ensure that it complies with the code rules of solvers like Gurobi, CPLEX, etc. 

\textbf{Step 3: Problem Tailor}. 
The main task of the problem tailor is to scope down the macro-level queries $Q$ into specific subsets of decisions $\hat{\textbf{x}}\in \mathcal{X}\backslash\textbf{x}^{\prime}$ and relevant factors rather than considering the entire collection of decisions and constraints in the E2E model $\min_{\textbf{x}\in \mathcal{X}} g(\textbf{x};\textbf{w})$.
The conversion is completed by using prompted $\textbf{LLM}(Q;\mathcal{D},P^{\text{PP}})$ based on empirical record analysis $\mathcal{D}$ and prompt $P^{\text{PP}}$. 
Specifically, as seen in Algorithm~\ref{pse:1}, the problem tailor first leverages LLM-generated SQL instructions to obtain compatible query-related information in databases $\mathcal{D}$, particularly those containing large-scale empirical operational data $(\bar{\textbf{x}},\bar{\textbf{w}})$. Furthermore, the problem tailor examines which variables are associated with the user's query based on historical data and subsequently generates the most relevant decision variables $\hat{\textbf{x}}$. 
This process can be formulated into an optimization problem to obtain the most relevant decision subset $\hat{\textbf{x}}$ about $Q$:
\begin{equation}
    \hat{\textbf{x}}_t=\arg\ \max_{\textbf{x}\in \mathcal{X}} \textbf{LLM}(Q;\mathcal{D},P_t^{\text{PP}}(\hat{\textbf{x}}_{t-1},S_{t-1})),
\end{equation}
where the prompt $P_t^{\text{PP}}(\hat{\textbf{x}}_{t-1},S_{t-1})$ in iteration $t$ involving previous iteration decisions $\hat{\textbf{x}}_{t-1}$ and satisfaction factor $S_t=\frac{f(\hat{\textbf{x}}^*_{t-1},\textbf{x}^{\prime};\textbf{w})-f(\hat{\textbf{x}};\textbf{w})}{f(\hat{\textbf{x}};\textbf{w})}$. 
Satisfaction factor $S_t$ is defined as the relative gap between the optimized QR-obj function value and the original value prior to optimization. It also enables LLM to identify the most relevant $\hat{\textbf{x}}$ and achieve a higher value of $S_{t}$.
For convenience of handling, $\textbf{x}^{\prime}$ are enforced to be equal to historical values $\bar{\textbf{x}}^{\prime}$ and are expressed as equality constraints in the following steps. 
After that, a code safeguard checks the LLM-generated constraints code to see if it complies with the solvers' syntax.

\textbf{Step 4: E2E Model}. The E2E program $\min_{\textbf{x}\in \mathcal{X}}g(\textbf{x},\textbf{w})$ (see details in the Section \ref{sec:e2e}) includes both decision variables for a specific problem and related features and parameters, called covariates. The E2E model includes two parts: (i) An offline trained Random Forest (Figure \ref{Fig.framework} \textcircled{4}) with a specific data predictor such as traffic congestion factors or operational cost of EBS. (ii) MIP: Following an earlier work \citep{Biggs_2022}, we encode the Random Forest tree structure as an MIP. The objectives of the MIP model include both operational objectives (identical with the RF data predictor) and query-related QR-objs $f(\textbf{x};\textbf{w})$ provided by QR-obj generators. 
In practice, the realization of query-relevant goal should not lead to unnecessary costs or compromise the original operational goals. Therefore, the optimization of the Query-relevant objective $f(\textbf{x};\textbf{w})$ will be tackled as a secondary priority, following the assurance of the optimality of the original objective in MIP. Additionally, the parameterized variables $\textbf{x}^{\prime}$, in the form of equality constraints, are incorporated into the MIP constraints. Then the agent proceeds with the following model to achieve scoped-down optimal solutions: 
\begin{equation}
\hat{\textbf{x}}^*_t=\arg \min_{\hat{\textbf{x}}} g(\hat{\textbf{x}};\bar{\textbf{x}}^{\prime},\textbf{w}).
\end{equation}
By parameterizing irrelevant decisions, the scoped-down problem $\min_{\hat{\textbf{x}}} g(\hat{\textbf{x}};\bar{\textbf{x}}^{\prime},\textbf{w})$ achieves efficient, optimal solutions $\hat{\textbf{x}}^*_t$ in less computational time. 

\textbf{Step 5: Response Prompter}.  
The E2E model iteratively generates optimal solutions until the user's query is satisfied ($S_t$). 
The satisfaction rate, which must be higher than zero and consistent across multiple iterations, controls the termination condition. Finally, the Response Prompter converts the optimal solutions and scores $(\textbf{x}^*,S^*)$ into a more easily comprehensible form for the end-users. 

\textbf{Pseudocode Summary}. 
We summarize the pseudocode of the City-LEO agents as follows: The pseudocode outlined in Algorithm~\ref{pse:1} encapsulates the operational flow of the City-LEO agent, designed to process users' queries input within an urban operations management context. 
The input comprises the natural language user query $Q$, an initial empty set $\textbf{x}^{\prime}_0$ of parameterized decisions, and a zero-initialized satisfaction score $S_0$, along with the LLM prompts $P^{\text{IG}}$ and $P^{\text{PP}}(\textbf{x}^{\prime}_0,S_0)$.
Initially, the Problem Matcher identifies an agent specific to the query's context. The QR-obj generator then starts, performing three critical tasks: it generates a query-relevant QR-obj function $f(\textbf{x};\textbf{w})$ via the LLM, a satisfaction factor $S_t$, and produces safeguarded objective function code.
The iterative process continues until either the maximum number of iterations is reached or the satisfaction factor becomes positive value (in the case of maximizing the QR-obj). 
Within each iteration, the problem tailor invokes empirical decision-related data  $\mathcal{D}$ from a database and updates the LLM-Optimal parameterized decisions $\textbf{x}^{\prime}_t$. Subsequently, the algorithm solves an $f$-objective embedded E2E programming  $\min_{\textbf{x}\in \mathcal{X}} g(\textbf{x};\textbf{w})$ that encapsulates the objective function $f$. It refines the solution space and updates both the satisfaction score $S_t$ and the LLM prompts $P^{\text{PP}}(\textbf{x}^{\prime}_0,S_0)$ for the next iteration.
Upon convergence or fulfillment of the termination criteria, the algorithm outputs the optimized decision vector $\textbf{x}^*$ alongside the final satisfaction score $S^*$, marking the end of the urban-centric City-LEO agent.
\begin{algorithm}[]
        \caption{City-LEO Agent Pseudocode} 
        \renewcommand{\algorithmicrequire}{\textbf{Input:}}
        \renewcommand{\algorithmicensure}{\textbf{Output:}}
        \renewcommand{\baselinestretch}{1.2}\selectfont
        \begin{algorithmic}
        \REQUIRE Query $Q$, initial value $S_0=0$, LLM prompts $P^{\text{IG}},P_0^{\text{PP}}(\textbf{x}^{\prime}_{0},S_{0})$.   
        \STATE \textbf{Problem Matcher:} Determine an area-specific agent.
        \STATE \textbf{QR-obj Generator:} 
        \STATE \hspace{10pt} 1. Query-relevant QR-obj function $f(\textbf{x};\textbf{w})\deq \textbf{LLM}(Q;P^{\text{IG}})$ generation.
                \STATE \hspace{10pt} 2. Objective $f(\textbf{x};\textbf{w})$ code generation (code safeguard checked).  
        \STATE \hspace{10pt} 3. Satisfaction factor $S_t=\frac{f(\hat{\textbf{x}}^*_t;\bar{\textbf{x}}^{\prime},\textbf{w})-f(\bar{\textbf{x}};\textbf{w})}{f(\bar{\textbf{x}};\textbf{w})}$ generation.
        \WHILE{$\{t\leq \text{max-iterations} \}$ or $\{S_t\textgreater 0\}$}
       \STATE \textbf{{a.}\ Problem Tailor:}
        \STATE \hspace{10pt} 1. Call decision-related empirical data $\mathcal{D} \deq (\bar{\textbf{x}},\bar{\textbf{w}})$ from a database.
        \STATE \hspace{10pt} 2. Update the LLM-Optimal parameterized decisions $\textbf{x}^{\prime}_t$: \\ 
        \begin{center}
            $\textbf{x}^{\prime}_{t+1}=\arg\ \max_{\textbf{x}\in \mathcal{X}} \textbf{LLM}(Q;\mathcal{D},P_{t}^{\text{PP}}(\textbf{x}^{\prime}_{t},S_{t})),$
        \end{center}
        \STATE \textbf{b.\ Solve a E2E model $g(\textbf{x};\textbf{w})$ embedded with QR-obj $f$: }
        \STATE \hspace{10pt} 1. Optimal Scoping down solutions: 
        \begin{center}
        $\hat{\textbf{x}}^*_{t+1}=\arg \min_{\hat{\textbf{x}}} g(\hat{\textbf{x}};\bar{\textbf{x}}^{\prime}_{t+1},\textbf{w}),$
        \end{center}
        \STATE \hspace{10pt} 2. Update $S_{t+1},P_{t+1}^{\text{PP}}(\textbf{x}^{\prime}_{t+1},S_{t+1})$.
        \STATE \hspace{10pt} \textbf{return} $\hat{\textbf{x}}^*_{t+1},S_{t+1}$.
        \ENDWHILE
        \RETURN $\left \{ \textbf{x}^* \right \} \left \{ S^* \right \}$.
	\end{algorithmic}\label{pse:1} 
\end{algorithm}

\subsection{Case: E-bike Sharing Operations}
To facilitate the understanding of the aforementioned agent framework, this section presents a practical example of how the proposed City-LEO can be applied in a real-world EBS system.
We specifically introduce prompts with the Q\&A instructions adopted from our City-LEO (as shown in Figure \ref{Fig.prompt1} and Figure \ref{Fig.prompt2}).
The query $Q$ pertains to operational-level pursuits commonly raised by EBS managers. 
We use $K=\{k_1, k_2, k_3\}$ to represent the set of discretized states of charge (SOC) levels of batteries, that is, $k_1, k_2$ and $k_3$ denote low SOC level (e.g. lower than 30\%), medium SOC level (e.g. 30\% \text{-} 70\%) and high SOC level (e.g. higher than 70\%), respectively.
Typical pursuits include reducing the number of low-powered bikes and adding available bikes at certain locations. 
In particular, Figure \ref{Fig.prompt1} shows the City-LEO agent working on a real-world operational problem and a user's question about EBS rebalancing.
This interaction involves both understanding and generating responses based on specific technical and domain-specific requirements. Next, we will elaborate on the user's query, prompt design, and thought-task breakdown as follows:

\textbf{User's Query}: The user, often an EBS manager, adopts City-LEO to inquire about operational decisions for increasing the number of available bikes at parking spots NO.1 and NO.2 (query $Q$) while ensuring the lowest possible operating cost. The crux of this matter lies in accurately defining "available bikes," leveraging operational data from the database in an automated manner, and identifying the most pertinent locations with NO.1 and NO.2 parking spots to re-optimize the EBS system. \\
\begingroup
\setlength\abovedisplayskip{0pt}
\setlength\belowdisplayskip{0pt}
\setlength\abovedisplayshortskip{0pt}
\setlength\belowdisplayshortskip{0pt}
\begin{figure}[]
\centering
\includegraphics[width=1\textwidth]{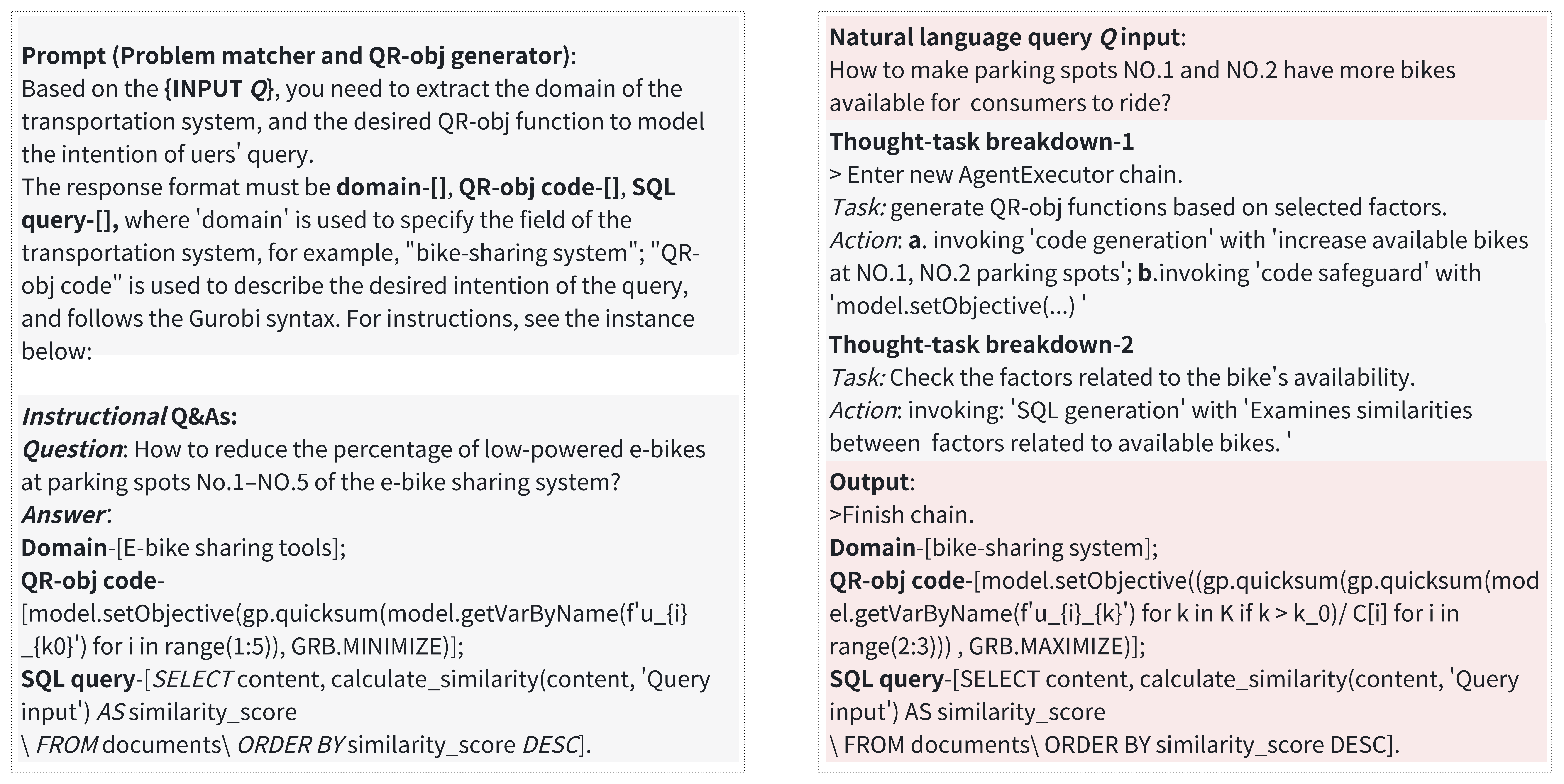}
\caption{Prompt in Problem Matcher and QR-obj Generator}
\label{Fig.prompt1}
\end{figure}
\endgroup
\textbf{Prompt (Problem Matcher and QR-obj Generator)}: The prompt first instructs LLM on how to extract the area-specific APIs (domain of the transportation system) and generate the desired QR-obj function from the user's query. The output response is specifically formatted as a "domain" and "QR-obj function codes" based on the prompt line (prompt matcher and QR-obj generators), as shown in Figures \ref{Fig.prompt1} and \ref{Fig.prompt2}. In addition, we also provide \emph{instructional Q\&A} prompts to improve the precision of generating SQL queries, domains, and codes of QR-obj functions.
In our case study, the given prompt specifies the domain as "EBS system," and the QR-obj code outputs a Gurobi solver code snippet. This code snippet can be utilized to optimize the number of available e-bikes at particular parking spots. 

\textbf{Thought Process Breakdown (Problem Matcher and QR-obj Generators)}: Based on the input and prompt, City-LEO breaks down the query into sub-tasks  following a step-by-step process: First, City-LEO analyzes the relevant factors associated with "available bikes" recorded in the database through LLM-generated SQL. 
These factors include historical records of operational-level decisions or parameters, such as rebalancing bikes and initial allocations. 
Secondly, City-LEO generates an QR-obj function and corresponding Gurobi code snippets (e.g., $\mathtt{model.setobjective}$(gp.quicksum(model.getVarByName(f'u\_\{i\}\_\{k\}') for i in range(55) for k in {[}1,2{]}), GRB.MAXIMIZE), see section \ref{sec:e2e} for notation details). This QR-obj function quantifies the number of "available bikes" using predetermined parameters in the initial stage. Subsequently, the code safeguard further guarantees the successful execution of the generated code in the E2E dispatch optimization model.


\textbf{Output}: The model's output consists of the domain "bike-sharing system" and codes for a QR-obj function aligned with Gurobi syntax. The QR-obj function, in particular, is meant to maximize the number of e-bikes with at least $k_1$ level state of charge (SOC) at No.1 and No.2 parking spots.

\textbf{Prompt (Problem Tailor):} After the area-specific agents have been matched, the users' inquiry is processed by the problem tailor.
The prompt of problem tailor consists of three components: historical data retrieval, problem scoping, and constraint generation. 
First, the prompt provides instructions on how to use SQL retrieve query-relevant operational information that is stored in the database. 
The information, such as the shared e-bike travel matrix and the redistribution of bikes among the spots, enables LLM to analyze the most relevant parking spots or decisions with queries. \\
\begingroup
\setlength\abovedisplayskip{0pt}
\setlength\belowdisplayskip{0pt}
\setlength\abovedisplayshortskip{0pt}
\setlength\belowdisplayshortskip{0pt}
\begin{figure}[]
      \centering
\includegraphics[width=1\textwidth]{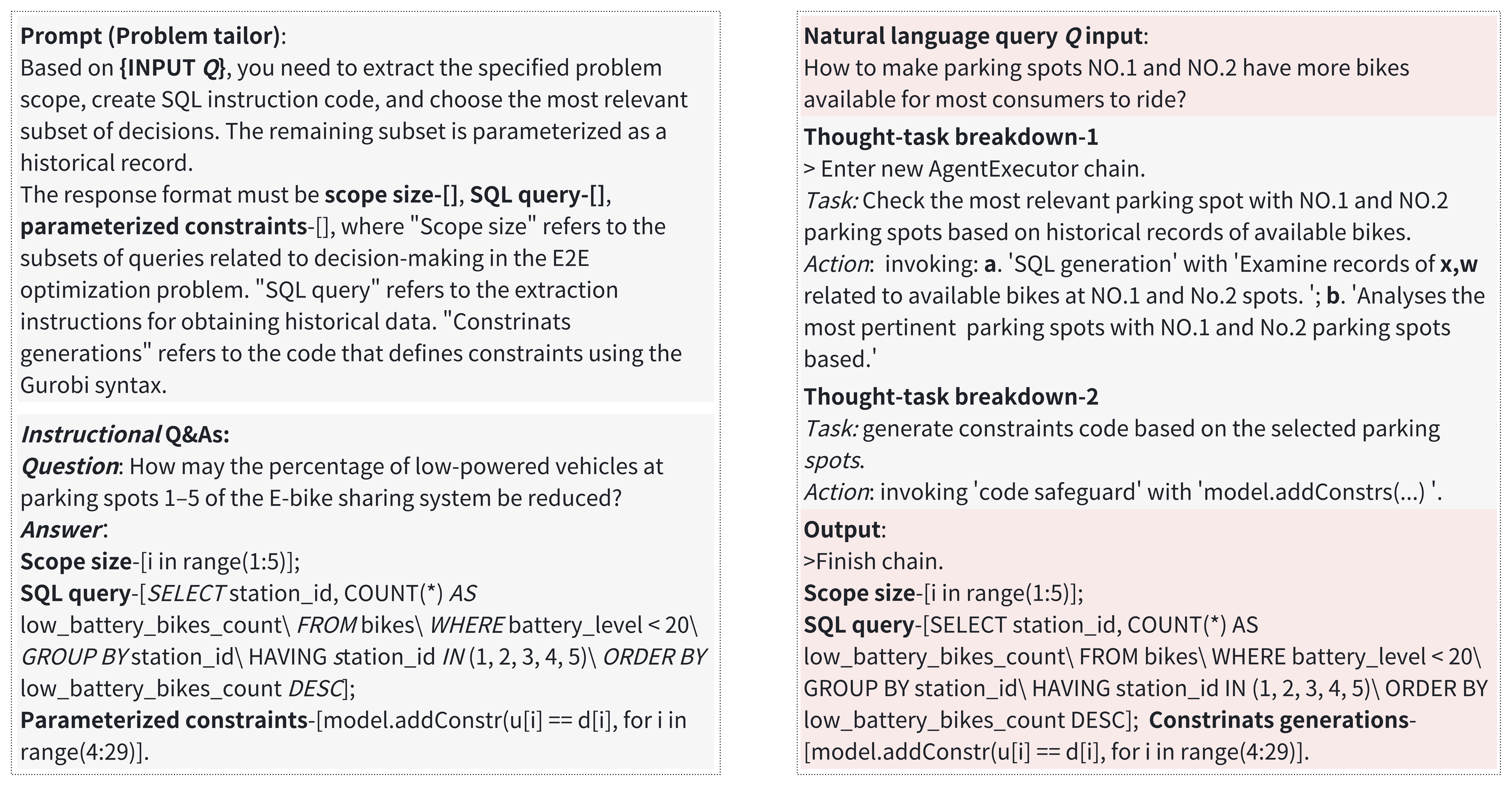}
\caption{Prompt in the Problem Tailor}
    \label{Fig.prompt2}
\end{figure}
\endgroup
For instance, Figure \ref{Fig.prompt2}  illustrates an example of an instructional Q\&A session on decreasing the proportion of low-powered vehicles at specified EBS parking spots.
Secondly, prompts instruct on problem scoping, specifically in this context, to aid in the selection of EBS parking spots that are relevant to the query. For example, parking spots ranging from No. 4 to No. 8 have a higher frequency of trips and a more frequent record of rebalancing operations from No. 1 and No. 2 parking spots. This helps to scope down the search area from $\mathcal{X}$ to a specific subset consisting of parking spots No. 1, No. 2, and No. 4 to No. 8 in the E2E dispatch optimization problem. 

\textbf{Thought Process Breakdown (Problem Tailor)}: Guided by the prompt of problem tailor, City-LEO then breaks down the user query into tasks including a query-relevant information search and constraint generation. In line with the task breakdown suggested in the QR-obj generator, the initial step involves retrieving relevant information, such as operational-level decisions or parameters, from a database using SQL. This information facilitates the process of determining the scope of a optimization problem and generating Gurobi constraints code (i.e., $\mathtt{model.addConstr}$()).

\textbf{Output:} The final output is a combination of SQL code, QR-obj and Gurobi optimization formulation tailored to the specifics mentioned in the query. This output is intended to facilitate decision-making that enhances bike availability based on historical usage patterns.
\section{End-to-End Model for the EBS System}\label{sec:e2e}
Having outlined the algorithm flow of City-LEO in Section \ref{sec:agent}, we next focus on the E2E model for the EBS system. EBS is an effective and appealing means of last-mile transportation by offering a convenient and flexible travel mode. 
However, the operators of the EBS system face challenges in dealing with the imbalance between supply and demand \citep{jin2023}, and recharging depleted batteries for e-bikes. 
To handle the operations management problems of urban EBS systems, in this section, we propose the optimization models embedded in the LLM agent that facilitate an efficient and transparent decision-making process. 
First, we construct the conventional deterministic optimization model for rebalancing e-bikes and battery swapping. Then we formulate the RF-based E2E optimization model, and leverage the feature-to-decision framework to obtain more accurate and interpretable solutions. In addition, to cope with diverse queries flexibly and reduce computational time, we explain how the LLM-based agent redefines the management problem and scopes down the size of the E2E optimization model.
\subsection{The E-bike Rebalancing and Battery Swapping Model}
Consider the EBS system with multiple parking spots, we construct an integer program to optimize the rebalancing and battery swapping strategies for the next several hours. At the initial of the horizon, two types of staff drive different vehicles to distribute e-bikes among multiple parking spots, and to replace depleted batteries with fully charged ones, respectively. Let $I=\{1,\cdots,|I|\}$ denote the set of parking spots. 
To ensure service quality, e-bikes with low SOC are not available for satisfying users' demands. Before rebalancing and battery switching, the initial number of e-bikes with SOC $k\;(k\in K)$ at parking spot $i$ is $v_{ik}$, and the capacity of parking spot $i$ is $C_i$. The capacity of the vehicle for carrying e-bikes is $N$, and the initial number of fully charged e-bikes on the dispatching vehicle is $r_{k_3}$. Other staff carries fully charged batteries to replace the depleted batteries, and the maximum number of fully charged batteries that the staff carries when setting off is $B$. 
After dispatching e-bikes and swapping batteries, the stock of e-bikes with SOC $k$ at each parking spot $i$ is $u_{ik}$. At the beginning, the operator has to decide the number of batteries $y_{ik}$ that are replaced by fully charged batteries at spot $i$, and the number of e-bikes $z_{ik}$ that are moved out or dropped off at spot $i$ with SOC $k$. In addition, let $c_1$ and $c_2$ denote the unit cost for swapping batteries and moving e-bikes, respectively. Define the customers' net demand at spot $i$ is $\xi_i$, and since customers always prefer to use e-bikes with high SOC, we add penalties $p_1$ for using e-bikes with medium SOC level to satisfy traveling demands. Penalties $p_2$ are for unmet demands $\sigma_i$. Given the notation, we formulate the following deterministic model to optimize the rebalancing and battery swapping (RBS) decisions:
\begin{align}
    \text{(RBS)}\quad \min\limits_{\mathbf{y,z,u}} \quad & \sum_{i} \left( \sum_{k<k_3} c_1 y_{ik} + \sum_{k\in K} c_2 |z_{ik}| \right) + p_1 \sum_{i} (\xi_i - \sigma_i - u_{ik_3})^+ + p_2 \sum_{i\in I} \sigma_i \label{obj1} \\
    \text{s.t.} \quad & u_{i k} = v_{i k} - y_{i k} + z_{i k}, \quad \forall i \in I, \; k < k_3, \label{con1}\\
    & u_{i k_3} = v_{i k_3} + \sum_{k < k_3} y_{i k} + z_{i k_3}, \quad \forall i \in I, \label{con2} \\
    & \sum_{k > k_1} u_{i k} + \sigma_i \geq \xi_i, \quad \forall i \in I, \label{con3} \\
    & \sum_{k \in K} \left( v_{i k} + z_{i k} \right) - \xi_i \leq C_i, \quad \forall i \in I, \label{con4}\\
    & \sum_{i \in I} \sum_{k < k_3} y_{i k} \leq B, \label{con5} \\
    & r_{k_3} - \sum_{i \in I} \sum_{k \in K} z_{i k} \leq N, \label{con6} \\
    & r_{k_3} - \sum_{i \in I} z_{i k_3} \geq 0, \label{con7} \\
    & \sum_{i} z_{i k_3} \geq 0, \label{con8} \\
    & \sum_{i} z_{i k} \leq 0, \quad \forall k < k_3, \label{con9} \\
    & u_{i k} \geq 0, \quad \forall i \in I, \forall k \in K; \; y_{ik} \geq 0, \quad \forall i \in I, \forall k \in K \backslash k_3. \label{con10}
\end{align}
Objective (\ref{obj1}) represents the total operational costs for swapping depleted batteries, dispatching e-bikes, and the penalties for unmet demands and not providing e-bikes with high SOCs to satisfy demands. Constraints (\ref{con1}) and (\ref{con2}) are the inventory-energy balance equations, that is, the stock of e-bikes with SOC $k$ at parking spot $i$ is equal to the stock in the previous period minus the number of replaced batteries and then plus the net dispatched e-bikes. Constraint (\ref{con3}) represents the relationship between the stock of e-bikes and the users' demands. Constraint (\ref{con4}) imposes space limits to make sure that the stock of e-bikes do not exceed the capacity of the parking spot. The initial number of fully charged batteries carried by the staff is $\sum_{i}\sum_{k<k_3}y_{ik}$, and constraint (\ref{con5}) sets the upper limit $B$ for the number of swapped batteries. Constraint (\ref{con6}) means that the number of e-bikes on the vehicle does not exceed the capacity limit $N$. Constraint (\ref{con7}) indicates that the total number of dispatched e-bikes is less than the initial number of fully charged e-bikes on the operator's vehicle. Constraint (\ref{con8}) indicates that the total number of fully charged e-bikes dispatched among the parking spots is nonnegative. Constraint (\ref{con9}) means that the total number of e-bikes with low SOC and medium SOC will not increase after rebalancing and battery swapping operations. Constraint (\ref{con10}) ensures the nonnegativity of the decision variables.
\subsection{RF-based E2E Model}
Considering the integration of OR models and LLM, the deterministic optimization problem (RBS) has some limitations. First, deterministic optimization cannot handle users' stochastic travel demands, which may result in suboptimal or even infeasible solutions. Second, the deterministic model leads to low data utilization and lacks the flexibility to efficiently respond to the diverse queries of users. 
Therefore, extending the approach of \citep{Biggs_2022}, we construct the following RF-based E2E model to involve more features and improve the utilization of historical data. In addition, the tree-based E2E optimization mitigates the error accumulation of predicting stochastic parameters, and avoids repeatedly training the RF models every time a user submits a query. The tree structure is also suitable for describing unknown and complex interactions of the covariates (e.g. features and variables) and objectives, and enhances the interpretability and transparency of the decision-making process.
\begingroup
\setlength\abovedisplayskip{0pt}
\setlength\belowdisplayskip{0pt}
\setlength\abovedisplayshortskip{0pt}
\setlength\belowdisplayshortskip{0pt}
\setlength\abovecaptionskip{0pt}
\setlength\belowcaptionskip{0pt}
\setlength\intextsep{5pt}
\begin{figure}[H]
     \centering
\includegraphics[width=0.9\textwidth]{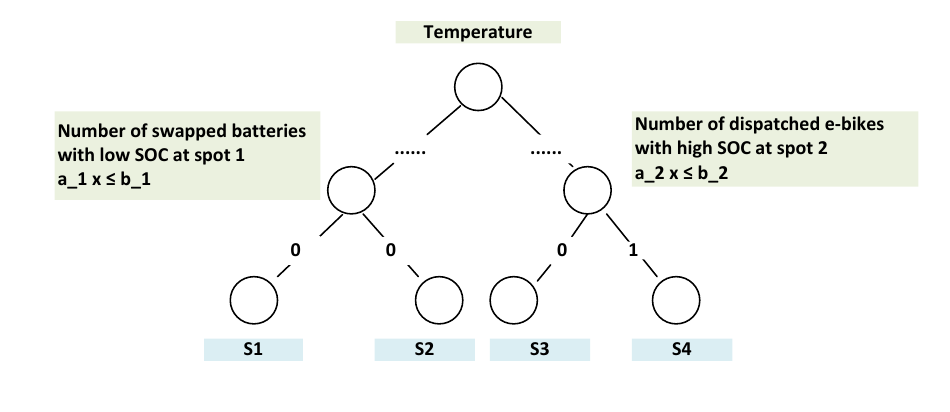}
   \caption{Sample Tree of the Trained Random Forest}
   \label{Fig.RF1}
\end{figure} 
\endgroup
For the EBS system, we use Figure \ref{Fig.RF1} to depict a sample tree of the RF, which is trained by the historical meteorological data and operational records. The tree structure can reflect the inherent uncertainty of users' demand. Suppose the number of trees in the RF is $|H|$, and the leaf node of each tree $h$ represents the predicted operational cost $R(\mathbf{x},\mathbf{w})$, where $\mathbf{x}=(\mathbf{y}, \mathbf{z})$ is the decision variable and $\mathbf{w}$ represents the features. Our goal is to minimize the cost for rebalancing and battery swapping by minimizing $R(\mathbf{x},\mathbf{w})$. Let $N^{h}$ denote the number of nodes (except the leaves) in tree $h$. For each interior node, let $p_i, l_i$ and $r_i$ be the immediate parent, the left child and the right child, respectively. We use $L^h$ to represent the
set of leaves in tree $h$, and use $R_j^{h}\;(j\in L^h)$ to denote the score of each leaf. Given the trained RF, we introduce binary variables $q_{ij}^{h}$ to select branches and decide the range of the variables, and the corresponding $0\text{-}1$ mixed integer programming is as follows:
\begin{align}
(\text{RBS-E2E})\;\min\limits_{\mathbf{x},\mathbf{q},\mathbf{u}}\quad & \frac{1}{H}\sum_{h=1}^{H}\sum_{j\in L^h}R_j^{h}q_{p_j,j}^{h} \label{obj2}\\
\text{s.t.}\quad & a_{i, h} \mathbf{x}-M\left(1-q_{i, l_i}^h\right) \leq b_i^h, \forall h \in H, i \in N^h,\label{con21} \\
& a_{i, h} \mathbf{x}+M\left(1-q_{i, r_i}^h\right) \geq b_i^h, \forall h \in H, i \in N^h,\label{con22}
\end{align}
\begin{align}
& q_{i, l_i}^h+q_{i, r_i}^h=q_{p_i, i}^h, \forall h \in H, i \in N^h,\label{con23} \\
& \sum_{i \in L^h} q_{p_i, i}^h=1, \forall h \in H,\label{con24}\\
& q_{i, l_i}^h, q_{i, r_i}^h, q_{p_i, i}^h \in\{0,1\}, \forall h \in H, i \in N^h, \label{con25} \\
& \text{Constraints} \; (\ref{con1})\text{-}(\ref{con2}), (\ref{con5})\text{-}(\ref{con10})\nonumber.
\end{align}
The objective (\ref{obj2}) is to minimize the predicted cost, that is minimize the score of the selected leaves. Constraint (\ref{con21}) and 
constraint (\ref{con22}) are the big-M reformulations of logical constraints and determine which leaf the solution $\mathbf{x}$ lies in. Constraint (\ref{con23}) ensures that if a parent node is inactive, its children must also be inactive; but if any child is active, then the parent must be active as well. Constraint (\ref{con24}) guarantees that within each tree $h$, one and only one leaf can be active. Constraint (\ref{con25}) defines the binary variables. According to the model (RBS), the EBS system also needs to satisfy constraints (\ref{con1})\text{-}(\ref{con2}) and (\ref{con5})\text{-}(\ref{con10}) to ensure capacity limitations, inventory-energy balance equations and so on.   

Compared with the deterministic optimization model (RBS), the model (RBS-E2E) integrates a random forest with MIP to involve more environmental features and uncertainties. 
The hierarchical tree structure can reflect the decision path from the root to a leaf node, which promotes the transparency and interpretability of the decision process. However, the model (RBS-E2E) lacks the flexibility to accommodate operators' diverse requirements in practical implementation. LLM offers a promising path to enhancing accessbility through conversational interactions. Therefore, we will introduce LLM-embedded E2E model in the subsequent section.


\subsection{LLM-embedded E2E Model}
The integration of the LLM agent and RF-based E2E optimization can lead to mutual enhancement for practical EBS management. According to the users' queries, the LLM agent utilizes its powerful reasoning ability and generates a new objective for the MIP. Users' queries usually focus on the operations of EBS system in various urban areas, and are less likely to involve the entire city. To cope with users' queries flexibly and avoid solving large-scale optimization models repeatedly, LLM contributes to scoping down the optimization problem to reduce computational complexity without significantly compromising the accuracy of the solutions. Specifically, after processing the users' queries, the LLM agent screens out the parking spots $\hat{I}$ and variables $\hat{\mathbf{x}}=(y_{ik}, z_{ik})\;(i\in \hat{I})$ that are most relevant to the query, and then parameterizes the other decision variables $\mathbf{x}^{\prime}=(y_{ik}, z_{ik})\;(i\in I\backslash \hat{I})$ with the mean values of the historical data $\bar{\mathbf{x}}^{\prime}$. For example, if a user is concerned about the operations of parking spots $\text{No}.5\text{-}\text{No}.10$, based on geographical distance and historical trip records, we can select some spots that are most relevant to the six spots (that is, the set of spots $\hat{I}$), rather than optimizing the EBS operations of the entire city. Figure \ref{Fig.RF2} represents that the variables associated with spot No.1 are not closely relevant to the user's problem. On the whole, the LLM-embedded optimization tool utilizes RF-based E2E optimization to learn from data and generate new decisions, and then provides responses to users through conversational dialogues. Define $A(\mathbf{w})F(\mathbf{x})\leq m(\mathbf{w})$ to represent the feasible region  $\mathbf{x}$. After scoping down the model (RBS\text{-}E2E), we have the following optimization program:
\begingroup
\setlength\abovedisplayskip{0pt}
\setlength\belowdisplayskip{0pt}
\setlength\abovedisplayshortskip{0pt}
\setlength\belowdisplayshortskip{0pt}
\begin{figure}[]
      \centering
\includegraphics[width=0.9\textwidth]{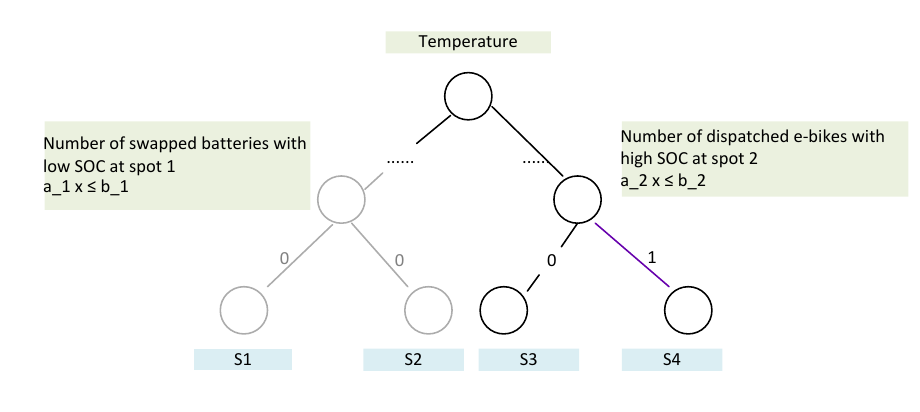}
    \caption{Sample Tree of the Trained Random Forest after Scoping Down the Problem}
    \label{Fig.RF2}
\end{figure} 
\endgroup
\begin{align}
(\text{RBS-LLM})\;\min\limits_{\hat{\mathbf{x}},\mathbf{q},\mathbf{u}} \quad & \frac{1}{H}\sum_{h=1}^{H}\sum_{j\in L^h}R_j^{h}q_{p_j,j}^{h} \label{obj3}\\
\min_{\hat{\mathbf{x}},\mathbf{q},\mathbf{u}}\quad &\textbf{LLM}(Q;\mathcal{D},P^{\text{PP}}) \label{obj4}\\
\text{s.t.}\quad & A(\mathbf{w})F(\hat{\mathbf{x}}, \bar{\mathbf{x}}^{\prime})\leq m(\mathbf{w}), \label{con31}\\
& \text{Constraints} \; (\ref{con23})\text{-}(\ref{con25})\nonumber.
\end{align}
Objective (\ref{obj4}) is generated by the LLM agent based on the user's query $Q$. To maintain minimum operational cost (\ref{obj3}) and then accommodate user's requirement,
we set a lower priority to objective (\ref{obj4}). In the multi-objective optimization problem, we solve the programming model with objective (\ref{obj3}) first, and then re-optimize the problem with the newly added objective (\ref{obj4}) while maintaining the first objective optimal. Constraint (\ref{con31}) represents the cluster of constraints about selecting branches of the trees, rebalancing e-bikes and battery swapping operations. Constraints (\ref{con23})\text{-}(\ref{con25}) involve binary variables to explain how to select branches and leaves in the random forest.

After scoping down the full-scale optimization problem and parameterizing some redundant variables $\textbf{x}^{\prime}$, the LLM-embedded E2E problem with shorter computational time is more flexible and transparent to cope with users' queries. 
In addition, the feasible region of model (RBS\text{-}LLM) gets smaller, and we will compare the optimality gaps and computational time in the case study.


\section{Case Study}\label{sec:CS}

In this section, we demonstrate the effectiveness of the City-LEO agent with a case study of the Pronto EBS system in Seattle, USA. 
Specifically, Section \ref{ssec:des} offers a brief description of the data and parameter settings. Section \ref{ssec:setup} introduces the experimental setup. Sections \ref{ssec:relev}, \ref{ssec:perfor}, and \ref{ssec:acc} present numerical results about relevance evaluation and accuracy evaluation of City-LEO, compared with Gurobi.
Furthermore, we assess the scope-down policy of City-LEO in the supplemental experiment outlined in Appendix \ref{app:scope}.

\subsection{Data Description and Parameter Settings}\label{ssec:des}
We adopt the real-world case of the Pronto Cycle Share system in Seattle. 
The system comprises 55 parking spots situated in the major neighborhoods of the city, with a fleet of 500 bicycles. 
Cycle Share Dataset \citep{cyclesharedataset} provided by Pronto contains data on the parking spots, trips, and weather information from October 13, 2014 to August 31, 2016. 
The trip records contain rental and return data, traveling distance, and users' information. The weather dataset encompasses daily weather-related indicators, such as temperature, humidity, and visibility. 
Based on the historical trip records and the stock information at parking spots, we generate data about e-bikes dispatching and battery swapping operations through simulation and ensure inventory-energy balance. 
Furthermore, sufficient historical data enhances City-LEO's reasoning ability for supplementing prior knowledge and generalization ability of E2E model. 

We analyze 283,143 pieces of historical trip records spanning from October 2014 to August 2016, and focus on the operations management of the EBS system during the noon hours, when the users are more active than the other periods. 
The capacities of vehicles for carrying e-bikes and batteries are 60 and 200, respectively. Assume that the initial number of e-bikes with high SOCs on the vehicle is $r_{k_3}=60$, and the initial number of fully charged batteries on aboard the vehicle is $B=142$. Regarding the operational costs, the unit cost for swapping depleted batteries or relocating e-bikes is \$1. In order to reduce the range anxiety and enhance the service quality, we set $p_1=0.4$ to penalize using e-bikes with medium SOC level to satisfy traveling demands, and the penalty parameter $p_2$ for stock-out of available e-bikes is \$0.8.
The cruising range of a high SOC level e-bike is approximately 50km-80km \citep{cyclereport}, and the average power consumption per trip is about 3\%. Therefore, SOC changes for each traveling trip is negligible. To validate the effectiveness of City-LEO across diverse datasets, we have further developed two additional datasets derived from Cycle Share Data. According to geographical proximity, we cluster the 55 parking spots into 20 and 35 clusters, respectively. Due to the page limit of the paper, we have omitted the detailed experimental results pertaining to the two datasets. The results consistently indicate that City-LEO has superior performance in terms of relevance and accuracy when benchmarked against the full-scale optimization model.




To train the random forest in the E2E optimization model, we import the historical trip data, geographic information of parking spots, weather information, and operational records into the LLM database, and set 80\% of the 8924 pieces of data as the training set, and the 20\% of the data is used for the test set. The maximum number of trees in a random forest is 100, and the maximum depth of trees is 400. The accuracy of random forest is 93.92\% on the training set and 58.22\% on the test set.


\subsection{Experimental Setup}\label{ssec:setup}
The goal of the case study is to assess City-LEO's ability to improve \textit{relevance} without excessively compromising \textit{accuracy} in resolving queries regarding EBS management. We conduct two types of numerical experiments, namely \textit{relevance} and \textit{accuracy} tests, in order to validate the advantages of City-LEO. 


In both experiments, we employ the unit-test approach, which is frequently utilized in software development, to assess the performance of City-LEO \citep{li2023large}.
We select three representative operational management scenarios within the EBS system, which comprehensively cover the perspectives of EBS operators, customers, and government.
For each scenario, we propose a set of 19 queries accompanied by human-labeled ground-truth answers (see Appendix \ref{app:gtqa}). 
To cope with the variability of the LLM's output, we conduct five experiments for each query and take the average values of the relevancy and accuracy assessments. 

We conduct relevancy and accuracy tests using GPT-4.0 and construct an agent framework based on Langchain (a typical framework for developing LLM-driven applications).
In addition, we solve the MIP in the E2E model using in Gurobi 10 on a macOS system with an Apple M1 Pro CPU and 16GB RAM. 


\subsection{Relevance Evaluation}\label{ssec:relev}
The objective function related to the user's query is proposed by utilizing LLM's reasoning ability on the input prior knowledge. To demonstrate the effectiveness of the generated objective function, we perform relevance tests to assess the similarity between the objective functions generated by City-LEO and ground-truth expressions labeled by human. 

In our comparative experiments, we not only scrutinize the expressions of the generated objectives, but also assess whether the generated function aligns with the ground-truth formulation in terms of mathematical essence. Sometimes the codes of the two functions are different, yet their mathematical implications remain identical. For example, the following two different codes both represent maximizing the total number of available e-bikes after dispatching and battery swapping operations.
\begin{itemize}
    \item \texttt{model.setObjective (gp.quicksum (model.getVarByName (f'u\_\{i\}\_\{k\}') for i in range(55) for k in [1, 2]), GRB.MAXIMIZE)}
    \item \texttt{model.setObjective (gp.quicksum (gp.quicksum (model.getVarByName (f'u\_\{i\}\_\{k\}') for k in K if k >= k\_2) for i in I), sense = GRB.MAXIMIZE)}
\end{itemize}

To this end, we propose two measures based on Jaro\text{-}Winkler distance, \textit{text similarity} and \textit{results similarity}. The two metrics evaluate the similarity between the "generated objective function" and the "ground-truth objective function" in terms of codes and mathematical essence (see Appendix \ref{app:essence} for details), respectively. Next, we analyze the performance of City-LEO for generating query-related objective functions under the two metrics.

\begingroup
\setlength\abovedisplayskip{0pt}
\setlength\belowdisplayskip{0pt}
\setlength\abovedisplayshortskip{0pt}
\setlength\belowdisplayshortskip{0pt}
\setlength\abovecaptionskip{0pt}
\setlength\belowcaptionskip{0pt}
\setlength\intextsep{5pt}
\begin{table}[H]
\centering
\caption{Results of the Relevance Test}\label{tab:rele}
 \renewcommand{\arraystretch}{0.7}
\begin{threeparttable} 
\begin{tabular}{@{}cccccc@{}}
\toprule
{Prompts\tnote{1}} & \multicolumn{2}{c}{{In-sample test}} & \multicolumn{2}{c}{{Out-of-sample test}} \\ \cmidrule(l){2-5} 
 & {Result similarity} & {Text similarity} & {Result similarity} & {Text similarity} \\ \midrule
0 & —— & —— & 0.61 & 0.90 \\
3 & 0.89 & 0.96 & 0.73 & 0.93 \\
6 & 0.91 & 0.96 & 0.78 & 0.94 \\
9 & 0.90 & 0.97 & —— & —— \\ \bottomrule
\end{tabular}
\begin{tablenotes}    
\footnotesize              
\item[1] The term "Prompts" refers to the number of Q\&A instructions used in each experiment. 
\end{tablenotes}          
\end{threeparttable}
\end{table}
\endgroup
Table \ref{tab:rele} illustrates that City-LEO performs well in out-of-sample evaluation test in which no instructions are given, with a similarity of 0.61 and 0.90 with the $0$ prompt.
Increasing the number of prompt sets (as shown in  "Prompts" in Table \ref{tab:rele}) mildly improves the similarity of the generated objective to the ground-truth function, both in the in-sample tests and the out-of-sample tests. The results demonstrate that City-LEO effectively leverage LLM's reasoning and generalization capabilities on prior knowledge to propose objective functions, which better match the users' queries and can be added to the E2E optimization model. With effective query-related objective, next we will investigate whether City-LEO has advantages in redefining and solving neighborhood-level mathematical programs, compared with the full optimization model for the whole city.
\subsection{Performance of City-LEO Operations}\label{ssec:perfor}
\begin{figure}[]
    \begin{minipage}[b]{0.5\textwidth}
      \centering
      \includegraphics[width=0.9\textwidth]{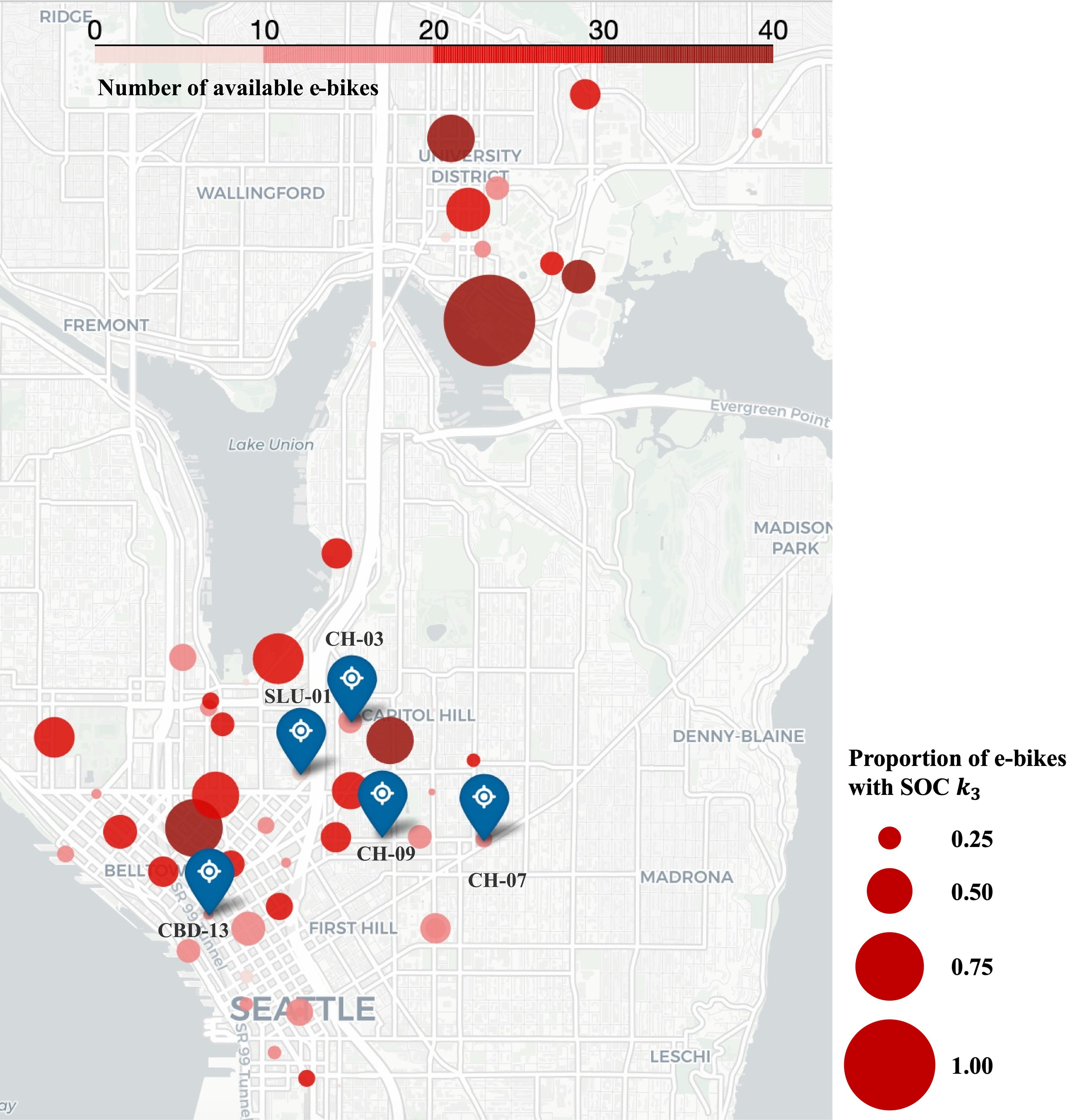}
      \vspace{-0.3cm}
    \begin{flushleft}
     \vspace{0.3cm}
      \small \centering (a) Distributional performance of FULL.
            \end{flushleft}
    \end{minipage}%
    \mbox{\hspace{0.12cm}}
    \begin{minipage}[b]{0.5\textwidth}
      \centering
      \includegraphics[width=0.9\textwidth]{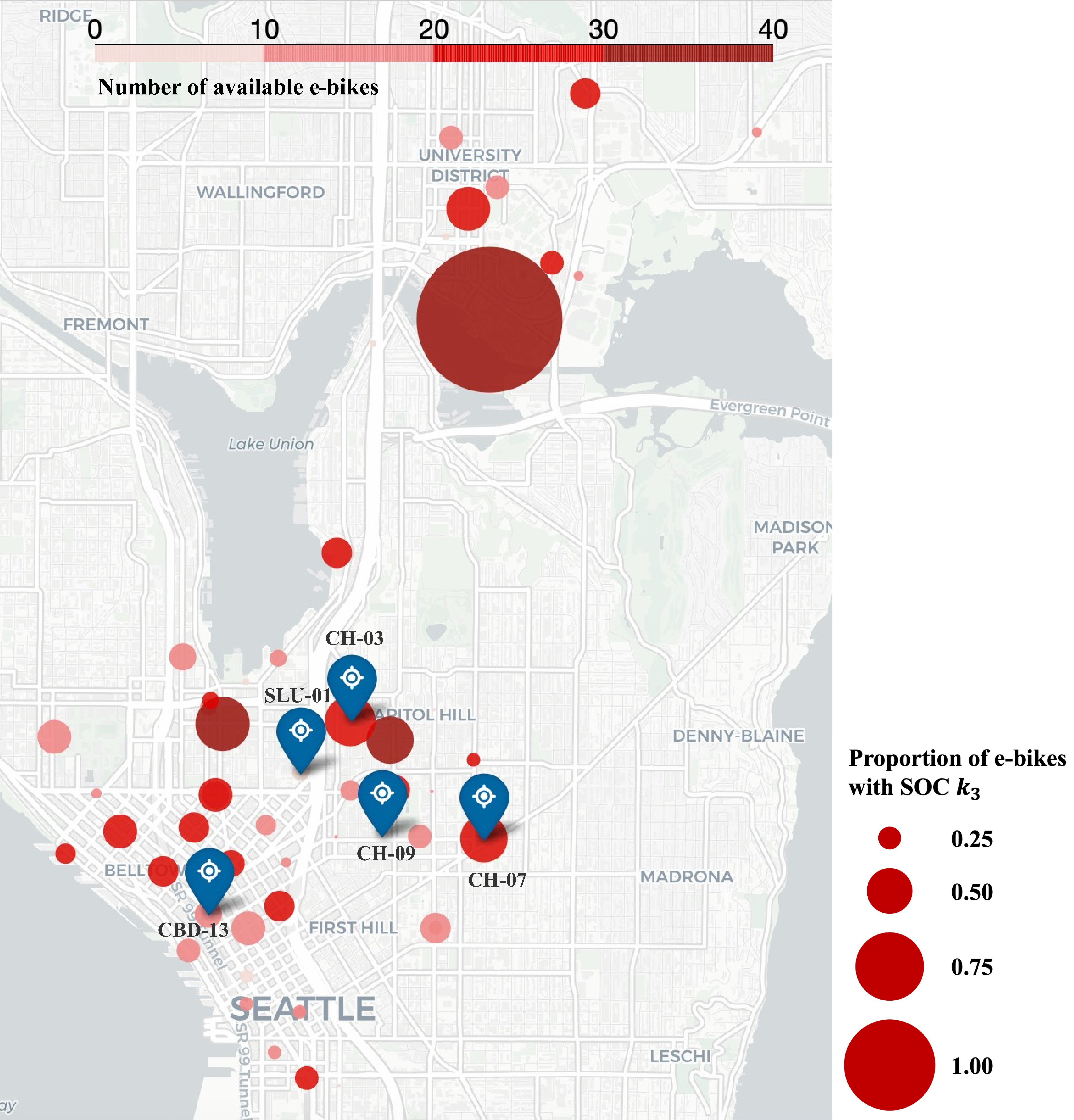}
      \vspace{-0.3cm}
     \begin{flushleft}
      \vspace{0.3cm}
      \small \centering (b) Distributional performance of City-LEO.
            \end{flushleft}
    \end{minipage}
    \caption{Spatial Distribution of EBS Operations Performance}
    \label{fig_comparison}
\end{figure}

In this section, we use City-LEO to enhance the quantity of available e-bikes, as well as increase the proportion of e-bikes with higher SOC, at the five designated parking locations shown by the blue symbols in Figure \ref{fig_comparison}.
We denote the baseline model as "FULL," representing the E2E optimization model without any scoping-down operations and aiming at optimizing the dispatch and battery swapping operations within the whole city. 
Figure \ref{fig_comparison} shows the spatial distribution of available e-bikes with SOC levels above $k_2$ and $k_3$ (size of dots), as well as the proportion of e-bikes at the $k_3$ SOC level (the varying shades of dots). More precisely, the bigger and darker shade of the dots indicates the presence of a greater number of available e-bikes and proportion of e-bikes at the $k_3$ SOC level.

We find that City-LEO outperforms FULL in satisfying users requirements. For example, as depicted in Figures \ref{fig_comparison}, for the five query-relevant parking spots, there are more available e-bikes with a higher SOC when adopting City-LEO (as shown in Figures \ref{fig_comparison}(b)) to optimize dispatch and battery swapping strategies.
The results highlight the practical value of employing City-LEO as an improved alternative to FULL in urban EBS operations.

\subsection{Accuracy Evaluation}\label{ssec:acc}
Compared with the FULL model, City-LEO aims to improve computational efficiency by inferring a scaled-down domain that is highly relevant to the user's query. For example, when a user submits a query regarding parking spots No.1-No.17 in an urban subregion, City-LEO needs to utilize the reasoning ability on the historical data to identify the parking spots that are most closely related to the operations of spots No.1-No.17 for optimization purposes. In the E2E optimization model, the decision variables associated with the less relevant parking spots are subsequently parameterized as the mean values of the historical data. As a result, our proposed method prioritizes the "scope" that is relevant to the users' queries, yet the optimization program offers suboptimal solutions compared to optimizing all the parking spots by the FULL model.
Hence, it is necessary to evaluate the discrepancies with the optimal solutions of the FULL model as well as computational efficiency.

City-LEO and the baseline FULL model share identical objectives, including the objectives regarding operational costs (referred to as "FULL-obj"), and the extra-generated query-relevant objectives (referred to as "QR-obj") in each experiment. 
We assess City-LEO's computational accuracy and effectiveness by examining two important deviations from the baseline model: "\textit{global suboptimality}" and "\textit{local satisfaction}," while also comparing CPU time.
The metrics of "global suboptimality" and "local satisfaction" quantify the differences between the optimal solution of the baseline mode "FULL" and City-LEO in terms of the "FULL-obj" and "QR-obj" values, respectively (see Appendix \ref{app:metric} for details). 

Furthermore, we propose the metric "Locality" to measure the proportion of parking spots mentioned in customers' queries. For instance, if customers only care about 17 specific parking spots among 55 EBS locations in total, then the locality is 30.91\%.
And the decision variables associated with these 17 designated parking spots and their relevant locations (analyzed by the City-LEO) will not be parameterized.
A higher level of locality indicates that the scoped-down problem is more similar to the FULL optimization model. 
In addition, we also conduct in-sample and out-of-sample tests to investigate the influence of prompts on the computational accuracy of City-LEO. 
Next, we will elaborate on our experimental results with query-relevant linear objectives and non-linear objectives, respectively.

\subsubsection{Performance of City-LEO with Linear Query-relevant Objective.}

As seen in Table \ref{tab:acc-LP}, compared with FULL model, City-LEO exhibits superior computational efficiency under different locality settings. In addition, City-LEO also has higher Local satisfaction values and possesses advantages in meeting users' requirements, while maintaining lower global suboptimality and not significantly compromising computational accuracy. For example, when the locality is 40\%, the user aims at optimizing the operations of 22 parking spots in a certain district of the city. Based on the geographic information and historical origin-destination information of the trips, City-LEO will screen out some parking spots that are operationally relevant to the 22 parking spots, and the variables associated with the other parking spots will be parameterized into historical values. 
The results show that the CPU time of City-LEO is about 16.26\% faster than FULL on average, the \emph{global suboptimality} is as small as 0.62\% and the \emph{local satisfaction gain} is 2.44\% on average. 

\begingroup
\setlength\abovedisplayskip{0pt}
\setlength\belowdisplayskip{0pt}
\setlength\abovedisplayshortskip{0pt}
\setlength\belowdisplayshortskip{0pt}
\setlength\abovecaptionskip{0pt}
\setlength\belowcaptionskip{0pt}
\setlength\intextsep{5pt}
\begin{table}[H]
\caption{Results of Accuracy Tests with Linear Query-relevant Objectives}\label{tab:acc-LP}
\centering
\resizebox{\textwidth}{!}{%
\begin{tabular}{@{}ccccccccc@{}}
\toprule
\multirow{4}{*}{Locality} & \multicolumn{4}{c}{In-sample test} & \multicolumn{4}{c}{Out-of-sample test} \\ \cmidrule(l){2-9} 
 & \multicolumn{2}{c}{CPU time (s)} & \multirow{3}{*}{\makecell{Global \\ suboptimality}} & \multirow{3}{*}{\makecell{Local \\ satisfaction gain}} & \multicolumn{2}{c}{CPU time (s)} & \multirow{3}{*}{\makecell{Global \\ suboptimality}} & \multirow{3}{*}{\makecell{Local \\ satisfaction gain} }\\ \cmidrule(lr){2-3} \cmidrule(lr){6-7}
 & \multirow{2}{*}{City-LEO} & \multirow{2}{*}{FULL} &  &  & \multirow{2}{*}{City-LEO} & \multirow{2}{*}{FULL} &  &  \\
 &  &  &  &  &  &  &  &  \\ \midrule
20\% & 168.56 & 195.64 & 0.71\% & 2.33\% & 154.96 & 179.17 & 0.98\% & 2.26\% \\
40\% & 171.65 & 212.57 & 0.62\% & 4.08\% & 183.86 & 218.90 & 0.65\% & 7.18\% \\
60\% & 180.97 & 212.83 & 0.59\% & 1.14\% & 179.91 & 209.45 & 0.72\% & 4.33\% \\
80\% & 178.34 & 214.83 & 0.56\% & 2.22\% & 180.59 & 199.84 & 0.95\% & 1.27\% \\ \bottomrule
\end{tabular}%
}
\end{table}
\endgroup
In addition, local satisfaction gain presents a comparison of "satisfaction levels" to users' requirements.  For instance, in Table \ref{tab:acc-LP}, the local satisfaction gains with varying locality values are all positive, which suggests that the EBS operation strategies provided by City-LEO better fulfill users' requirements compared to FULL. 
Furthermore, in the out-of-sample tests, the \emph{global suboptimality} and \emph{local satisfaction gain} are consistently maintained at lower levels, which suggests that City-LEO is less prone to overfitting when faced with various queries.

\subsubsection{Performance on Nonlinear Query-relevant Objectives.}
City-LEO sometimes generates convex nonlinear objective functions based on the intent of users' queries. For example, City-LEO generates ${\sum_{i\in I}(\sum_{k\in K}u_{ik}}-{C_i})^2$  to denote the \textit{utilization rate of parking spots} or $\sum_{i\in I }\sum_{k\in \{k_1,k_2\}} (x_{ik})^2$  to present the turnover rate of e-bikes with a low SOC.
Embedding such a quadratic objective into the E2E optimization model often results in an inefficient computational process.
In order to assess the performance of City-LEO with nonlinear objectives in out-of-sample tests, we offer two representative instances, namely "maximizing the \textit{utilization rate of parking spots}" (denoted as "NLP-obj-1") and "minimizing the \textit{turnover rate of e-bikes with low SOC} (denoted as "NLP-obj-2"), respectively.

Considering two different nonlinear objectives, Table \ref{tab:acc-NLP} shows the comparison results of computational performances between City-LEO and FULL. Results show that City-LEO exhibits superior computational efficiency while incurring certain compromises on accuracy. 
For instance, with generated objective NLP-obj-1, City-LEO takes only 94.92s to solve the E2E optimization problem, while the FULL model fails to generate the optimal solutions within one hour (and the optimality gap is about 66\%). 
We compare the optimal solutions provided by City-LEO with the feasible solutions generated by FULL after one-hour computation. When the user is concerned about 20\% of the parking spots, the global suboptimality is 11.75\%, and the local satisfaction gain -3.70\%.  
\begingroup
\setlength\abovedisplayskip{0pt}
\setlength\belowdisplayskip{0pt}
\setlength\abovedisplayshortskip{0pt}
\setlength\belowdisplayshortskip{0pt}
\setlength\abovecaptionskip{0pt}
\setlength\belowcaptionskip{0pt}
\setlength\intextsep{5pt} 
\setlength\dblfloatsep{0pt}
\begin{table}[H]
\caption{Results of Accuracy Tests with Nonlinear Query-relevant Objectives}\label{tab:acc-NLP}
\centering
\resizebox{\textwidth}{!}{%
\begin{tabular}{@{}ccccccccc@{}}
\toprule
\multirow{4}{*}{Locality} & \multicolumn{4}{c}{NLP-obj-1} & \multicolumn{4}{c}{NLP-obj-2} \\ \cmidrule(l){2-9} 
 & \multicolumn{2}{c}{CPU time(s)/GAP(\%)} & \multirow{3}{*}{\makecell{Global \\ suboptimality}} & \multirow{3}{*}{\makecell{Local \\ satisfaction gain}} & \multicolumn{2}{c}{CPU time(s)} & \multirow{3}{*}{\makecell{Global \\ suboptimality}} & \multirow{3}{*}{\makecell{Local \\ satisfaction gain}} \\ \cmidrule(lr){2-3} \cmidrule(lr){6-7}
 & \multirow{2}{*}{City-LEO} & \multirow{2}{*}{FULL} &  &  & \multirow{2}{*}{City-LEO} & \multirow{2}{*}{FULL} &  &  \\
 &  &  &  &  &  &  &  &  \\ \midrule
20\% & 94.92 & 3600 (66.44\%) & 11.75\% & -3.70\% & 85.21 & 520.99 & 11.96\% & 6.65\% \\
40\% & 149.19 & 3600 (66.76\%) & 8.14\% & 0.07\% & 90.24 & 485.88 & 7.86\% & 4.88\% \\
60\% & 468.64 & 3600 (66.16\%) & 5.23\% & -4.35\% & 158.81 & 499.04 & 6.67\% & -7.50\% \\
80\% & 624.80 & 3600 (66.37\%) & 2.08\% & -0.92\% & 315.18 & 486.25 & 4.03\% & -4.32\% \\ \bottomrule
\end{tabular}
}
\end{table}
\endgroup
In terms of NLP-obj-2, we observe that City-LEO exhibits similar benefits in terms of CPU time. As the level of localization increases, City-LEO's Global suboptimality progressively decreases, e.g., from 11.96\% to 4.03\%. When the locality is 20\% or 40\%, the results indicate that City-LEO better satisfies the user's requirements. 
In summary, when implementing QR-obj with quadratic terms, City-LEO significantly outperforms the baseline FULL.

\section{Conclusion}\label{sec:con}
In this paper, we introduce City-LEO, an LLM agent integrated with E2E optimization for enhancing practical relevance and decision-making efficacy in an urban EBS management case. 
Such an agent is designed to bridge the gap between the intricate nature of OR tools and the practical needs of city operators (EBS operators in our context). It provides a conversational interface that realizes the effective accomplishment of query-relevant quantitative targets by flexibly scoping down large-scale optimization problems. It also offers responsible tools which are more user-friendly, transparent for city operators with little knowledge of operations research.
Furthermore, by leveraging pre-trained E2E optimization models and incorporating stakeholder natural language query-relevant features, City-LEO also enhances the quality of solutions in uncertain urban decision environments.  

In a real-world case, we examine the performance of our proposed City-LEO in terms of query relevance and optimality for large-scale problems. 
In particular, we observe that City-LEO shows reliable query\text{-}relevant objective generation, achieving 0.90 text similarity without prompts in out-of-sample test. The result ensures that subsequent E2E models obtain optimal solutions closely relevant to the user's query.
In the accuracy test, City-LEO demonstrates better computational efficiency, particularly when dealing with nonlinear query-relevant objectives, and achieves lower levels of global suboptimality relative to the baseline model under various types of query-relevant objectives. 
The results indicate that our proposed City-LEO achieves greater relevance with users' requirements while incurring reasonably small compromise on optimality.
At last, through adjusting the maximum permitted number of parameterized parking spots, City-LEO has been demonstrated to provide a reasonable range of parameterization operations to fulfill users' queries.
The aforementioned experimental advantages essentially stem from the benefits of \emph{scoping down} operations inside E2E model via City-LEO. LLM's reasoning capability enables City-LEO to make substantial progresses in reducing the complexity and intransparency, which are typically encountered in traditional OR tools used for urban operations management. 

In addition, by facilitating a clearer understanding of the decision-making process, City-LEO helps align municipal decisions with community needs and enhances public trust in an interactive way.
Furthermore, the outcomes of our study not only affirm the viability of using LLMs and E2E optimization in public administration but also reveal the potential for the technologies to transform city management into a more inclusive and transparent practice. 
While City-LEO significantly advances the state of the art, future work may explore integrating more dynamic data inputs and real-time decision-making capabilities to further enhance its responsiveness and accuracy. Additionally, extending the framework to incorporate feedback loops from implemented decisions could refine its prescriptive accuracy and user satisfaction, creating a continuously improving system that adapts to the evolving urban landscape.

\bibliographystyle{informs2014}  
\bibliography{literatureall}  

\clearpage
\appendix
\section{Notation}\label{app:notation}
For the readers' convenience, we provide the notations used in City-LEO agent and E2E models as follows:
 \renewcommand{\thetable}{4}
\begin{table}[H]
\centering
\caption{Notations Table}\label{tab-app:notation}
\resizebox{\textwidth}{!}{%
 \renewcommand{\arraystretch}{1.2}
\begin{tabular}{l l}
\hline
\textbf{Symbol Description} & \\ \hline
$Q$ &  Natural language query proposed by users.\\ 
$f(\textbf{x};\textbf{w})$ &  QR-obj function. \\ 
$g(\textbf{x};\textbf{w})$ &  E2E programming. \\ 
$\bar{\textbf{x}}$ & Historical record of decisions. \\ 
$\bar{\textbf{w}}$ & Historical record of parameters.\\ 
$\mathcal{D}$ & Empirical data $\mathcal{D}\deq(\bar{\textbf{x}},\bar{\textbf{w}})$ recorded in Database.\\ 
$P^{\text{IG}}$ & Prompts of QR-obj generator with Q\&A instructions.\\ 
$P^{\text{PP}}$ & Prompts of problem tailor with Q\&A instructions.\\ 
$\textbf{x}^{\prime}$ & Parameterized decisions. \\ 
$\hat{\textbf{x}}$ & Non-parameterized decisions. \\ 
$S_t$ & Satisfaction factor in $t$-th iteration. \\ \hline
 \textbf{Sets} & \\ \hline
 $\mathcal{X}$ &  Feasible domain of decision variables.\\
$\mathcal{W}$ & Parameters set.\\ 
$I$  & Set of parking spots.\\
$K$  & Set of discretized SOC levels of e-bike batteries.\\\hline
\textbf{Parameters} & \\ \hline
$\textbf{w}$ &  Parameters $\textbf{w}\in \mathcal{W}$. \\ 
$N$  & Capacity of the vehicle for dispatching e-bikes.\\
$B$  & Capacity of the vehicle carrying fully-charged batteries.\\
$C_i$  & Capacity of parking spot $i$.\\
$v_{ik}$ & Initial number of e-bikes with SOC $k$ at parking spot $i$.\\
$r_{k_3}$ &  Initial number of fully charged e-bikes on the dispatching vehicle.\\
$\xi_i$ & Net traveling demand at parking spot $i$.\\
$\sigma_i$  & Unmet traveling demand at parking spot $i$.\\
$p_1$ & Penalty for using e-bikes with medium SOC level to satisfy demands.\\
$p_2$ & Penalty for unmet demands.\\
\hline
\textbf{Decision Variables}&\\\hline
$\textbf{x}$ &  Decisions variables $\textbf{x}\in \mathcal{X}$. \\ 
$u_{ik}$ & Stock of e-bikes with SOC $k$ at parking spot $i$ after dispatching and battery swapping. \\
$y_{ik}$ & Number of batteries with SOC $k$ that are replaced by fully charged batteries at parking spot $i$.\\
$z_{ik}$ & Number of e-bikes that are moved out or dropped off at parking spot $i$.\\
\hline
\end{tabular}
 \renewcommand{\arraystretch}{1}
 }
\end{table}

\section{Experiments Settings}
\subsection{Ground-truth Q\&A}\label{app:gtqa}
In this section, we provide the ground-truth Q\&A instructions adopted in the prompt of City-LEO. 
Specifically, we provide accurate Q\&A information (human-labeled) in three specific scenarios: \emph{Operational}, \emph{Customer-related}, and \emph{Regulatory}. These scenarios encompass the majority of operational inquiries proposed by EBS operators. 
Within each scenario, there are typically 3 to 8 key QR-objs that are of utmost relevance to the majority of EBS operators (as shown in Table. \ref{app-tab:groundtruth}). 
For each QR-obj, we provide a ground-truth "Answers" that is derived from human experience. This information is provided in the form of a Gurobi objective code.

 \renewcommand{\thetable}{5}
\begin{sidewaystable}
  \caption{Ground-truth Objectives}\label{app-tab:groundtruth}
  \vspace{0.2cm}
  \centering
  
  \renewcommand{\arraystretch}{2}
  \begin{adjustbox}{max width=\textwidth,center}
  \resizebox{\linewidth}{!}{%
  \begin{tabular}{c c l}
    \hline
    \textbf{Scenario} & \textbf{QR-obj} & \textbf{Ground-truth Obj-code} \\ \hline
    \multirow{7}{*}{\textbf{Operations}} & \makecell{Proportion of bikes \\ with lower SOC} & $\mathtt{model.setObjective(gp.quicksum(model.getVarByName(f'u\_\{i\}\_\{k0\}')\ for\ i\ in\ range(55)), GRB.MINIMIZE)}$ \\
    & Utilization of bike & $\mathtt{model.setObjective(gp.quicksum(model.getVarByName(f'u\_\{i\}\_\{k\}')   for\ i\ in\ range(55)\ for\ k\ in\ {[}1,\ 2{]}), GRB.MAXIMIZE)}$ \\
    & Charging cost & $\mathtt{model.setObjective(gp.quicksum(model.getVarByName(f'x\_\{i\}\_exchange\_\{k\}')   * (2 - k)\ for\ i\ in\ range(55)\ for\ k\ in\ {[}0,\ 1{]}), GRB.MINIMIZE)}$ \\
    & \makecell{Daily average battery \\ exchange frequency} & $\mathtt{model.setObjective(gp.quicksum(model.getVarByName(f'x\_\{i\}\_exchange\_\{k\}')\ for\ i\ in\ range(55)\ for\ k\ in\ {[}0,\ 1{]}), GRB.MINIMIZE)}$ \\
    & Market share & $\mathtt{model.setObjective(gp.quicksum(model.getVarByName(f'u\_\{i\}\_\{k\}')\ for\ i\ in\ range(55)\ for\ k\ in\ {[}1,\ 2{]}), GRB.MAXIMIZE)}$ \\
    & Order cancellation rate & $\mathtt{model.setObjective(gp.quicksum(model.getVarByName(f'u\_\{i\}\_\{k\}')\ for\ i\ in\ range(55)\ for\ k\ in\ {[}1,\ 2{]}), GRB.MAXIMIZE)}$ \\
    & Market penetration rate & $\mathtt{model.setObjective(gp.quicksum(model.getVarByName(f'u\_\{i\}\_\{k\}')\ for\ i\ in\ range(55)\ for\ k\ in\ {[}1,\ 2{]}), GRB.MAXIMIZE)}$ \\ \hline
    \multirow{4}{*}{\textbf{Regulatory}} & overparking rate & $\mathtt{model.setObjective(gp.quicksum(gp.quicksum(model.getVarByName(f'u\_\{i\}\_\{k\}')\ for\ k\ in\ K)\ /\ C{[}i{]}\ for\ i\ in\ I)\ /\ len(I), GRB.MINIMIZE)}$ \\
    & Carbon emission reduction & $\mathtt{model.setObjective(gp.quicksum((2-k)   * model.getVarByName(f'x\_\{i\}\_exchange\_\{k\}') for\ i\ in\ range(55)\ for\ k\ in\ {[}0,\ 1{]}), GRB.MAXIMIZE)}$ \\
    & Customer adoption & $\mathtt{model.setObjective(gp.quicksum(model.getVarByName(f'u\_\{i\}\_\{k\}')\ for\ i\ in\ range(55)\ for\ k\ in\ {[}1, 2{]}), GRB.MAXIMIZE)}$ \\
    & Service coverage area & $\mathtt{model.setObjective(gp.quicksum(model.getVarByName(f'u\_\{i\}\_\{k\}')\ for\ i\ in\ range(55)\ for\ k\ in\ {[}1, 2{]}), GRB.MAXIMIZE)}$ \\ \hline
    \multirow{8}{*}{\textbf{Customer}} & Complaint rate & $\mathtt{model.setObjective(gp.quicksum(k*gp.quicksum(model.getVarByName(f'u\_\{i\}\_\{k\}')\    for\ i\ in\ range(55))\ for\ k\ in\ {[}1, 2{]}), GRB.MAXIMIZE)}$ \\
    & Quality of service & $\mathtt{model.setObjective(gp.quicksum(model.getVarByName(f'u\_\{i\}\_\{k\}')\ for\ i\ in\ range(55)\ for\ k\ in\ {[}1, 2{]}), GRB.MAXIMIZE)}$ \\
    & Customers' satisfaction level & $\mathtt{model.setObjective(gp.quicksum((k+1)   * model.getVarByName(f'u\_\{i\}\_\{k\}')\ for\ i\ in\ range(55)\ for\ k\ in\ {[}1, 2{]}),   GRB.MAXIMIZE)}$ \\
    & \makecell{Accessible bike\\ with higher SOC} & $\mathtt{model.setObjective(gp.quicksum(model.getVarByName(f'u\_\{i\}\_\{k\}')\ for\ i\ in\ range(55)\ for\ k\ in\ {[}1, 2{]}), GRB.MAXIMIZE)}$ \\
    & Average riding distance & $\mathtt{model.setObjective(gp.quicksum(model.getVarByName(f'u\_\{i\}\_\{k\}')*k\ for\ i\ in\ range(55)\ for\ k\ in\ {[}0, 1, 2{]}), GRB.MAXIMIZE)}$ \\
    & Customers retention rate & $\mathtt{model.setObjective(gp.quicksum(model.getVarByName(f'u\_\{i\}\_\{k\}')\ for\ i\ in\ range(55)\ for\ k\ in\ {[}0, 1, 2{]}), GRB.MAXIMIZE)}$ \\
    & Under growth rate & $\mathtt{model.setObjective(gp.quicksum(model.getVarByName(f'u\_\{i\}\_\{k\}')\ for\ i\ in\ range(55)\ for\ k\ in\ {[}1, 2{]}), GRB.MAXIMIZE)}$ \\
    & \makecell{Proportion of \\ available bike} & $\mathtt{model.setObjective((gp.quicksum(gp.quicksum(model.getVarByName(f'u\_\{i\}\_\{k\}')\  for\ k\ in\ K\ if\ k\ \textgreater\ k\_0)\ /\  C{[}i{]}\ for\ i\ in\ I)) , GRB.MAXIMIZE)}$ \\ \hline
  \end{tabular}
  }
  \end{adjustbox}
\end{sidewaystable}
\subsection{Sample of Gurobi Code Created by City-LEO}\label{app:essence}
This section presents the mathematical representation of the Gurobi code created by City-LEO, which is utilized to quantify the "Results Similarity" in the relevance test.
The Gurobi code is shown below, which represents the objective of maximizing the total number of accessible e-bikes and its related mathematical essence:

\begin{itemize}
    \item \textbf{Gurobi code}: \texttt{model.setObjective(gp.quicksum(model.getVarByName(f'u\_\{i\}\_\{k\}') for i in range(55) for k in [k\_2, k\_3]), GRB.MAXIMIZE)}
    \item  \textbf{Mathematical essence}:  $\max\  u_{0,1}+u_{0,2}+u_{1,1}+\underbrace{\quad\quad\quad \quad\cdots\cdots\quad\quad\quad\quad}_{u_{i,k}\ for\  i\ in\ [2,53]\ for\ k\ in\ [1,2]}+u_{54,1}+u_{54,2}$
    \item \textbf{Gurobi code}: \texttt{model.setObjective(gp.quicksum(gp.quicksum(model.getVarByName( f'u\_\{i\}\_\{k\}') for k in K if k >= k\_2) for i in I), sense = GRB.MAXIMIZE)}
    \item  \textbf{Mathematical essence}: $\max\ u_{0,1}+u_{0,2}+u_{1,1}+\underbrace{\quad\quad\quad \quad\cdots\cdots\quad\quad\quad\quad}_{u_{i,k}\ for\  i\ in\ [2,53]\ for\ k\ in\ [1,2]}+u_{54,1}+u_{54,2}$
\end{itemize}

The aforementioned two scenarios demonstrate that although City-LEO generates different Gurobi codes, the two codes nevertheless possess the same underlying mathematical essence.


\subsection{Metrics of Accuracy Tests}\label{app:metric}
We propose two criteria, \textit{Global suboptimality} and \textit{Local satisfaction gain}, to assess City-LEO's performance in terms of accuracy.
Recall that FULL-obj refers to the initial objective function of the baseline model (i.e., FULL). Also, it represents the aggregate value of the leaf nodes in the random forest within the E2E model. The FULL-obj additionally represents the basis of EBS operations, with the goal of minimizing operating costs. Yet, achieving the second purpose (referred to as QR-obj) may conflict with the FULL-obj. Thus, we can adopt \textit{global suboptimality} to analyze the compromising degree of FULL-obj when achieving users' QR-obj. 
Specifically, \textit{global suboptimality} ($g^{\text{s}}$) can be formulated as follows:

\begin{align*}
    g^{\text{s}}\deq\frac{\text{FULL-obj}_{\text{City-LEO}}-\text{FULL-obj}_{\text{FULL}}}{\text{FULL-obj}_{\text{FULL}}}
\end{align*}
where $\text{FULL-obj}_{\text{City-LEO}}$ and $\text{FULL-obj}_{\text{FULL}}$ represent the optimal objectives of City-LEO and FULL in terms of initial objective function, respectively. A lower value of $g^{\text{s}}$ suggests a lower level of compromise in attaining consumers' QR-obj with FULL-obj.

Secondly, the \textit{local satisfaction gain} quantifies the "gains of QR-obj" achieved by the City-LEO model in terms of parking spots (refereed as the word "local") mentioned in user inquiries compared to the FULL. The FULL model does not involve any scoped-down operation. This measure ($g^{\text{l}}$) is formulated as follows:
\begin{align*}
    g^{\text{l}}\deq\frac{\text{QR-obj}_{\text{City-LEO}}-\text{QR-obj}_{\text{FULL}}}{\text{QR-obj}_{\text{FULL}}}
\end{align*}
where $\text{QR-obj}_{\text{City-LEO}}$ and $\text{QR-obj}_{\text{FULL}}$ represent the optimal values of QR-obj generated by City-LEO and FULL, respectively. These values are calculated based on the parking spots specified in user inquiries. 
Moreover, QR-obj's optimization direction (maximum or minimum) is frequently linked to the user's inquiry. For instance, increasing the number of available bikes while decreasing complaints represents a distinct path. In this regard, we separately computed the value of $\text{QR-obj}_{\text{City-LEO}}-\text{QR-obj}_{\text{FULL}}$ based on different QR-obj types. Therefore, if $g^{\text{l}}\geq 0$, it means City-LEO outperforms FULL, and vice versa.

\section{Sample of Prompts}
This section provides detailed information regarding prompt settings in the City-LEO agent architecture (See, Figure \ref{Fig.framework}), which are categorized into Problem matcher, QR-obj generators, and Problem tailor.
\subsection{Problem Matcher}
The purpose of the Problem matcher is to match precise APIs that are encapsulated with corresponding QR-obj generators and E2E models. Therefore, the prompt only emphasizes instructing City-LEO to locate the most query-relevant APIs. The corresponding prompts are listed as follows:\\
\textbf{Function description:}
\begin{flushleft}
\setlength{\rightskip}{0pt plus 1fil}
\raggedright
\texttt{
Extract the transportation domain and target from user's input.
}
\end{flushleft}
\textbf{Prompt template:}
\begin{flushleft}
\setlength{\rightskip}{0pt plus 1fil}
\raggedright
\texttt{
User's questions are as follows:\\
\textit{\{Input\}}\\
Based on the user's input, you need to extract the domain and the target that the user is concerned about. It is necessary to ensure that the output adheres to the prescribed format. Pay close attention to the details of the user's target, and extract the information directly without summarizing it.
}
\end{flushleft}
\subsection{QR-obj Generator}
The purpose of the QR-obj generators is to formulate the user's target into a objective function that is convex with decision variables and parameters, which are defined in the OR model. The formula is then translated into Gurobi syntax. The corresponding prompts are listed as follows:\\
\textbf{Function description:}
\begin{flushleft}
\setlength{\rightskip}{0pt plus 1fil}
\raggedright
\texttt{
Generate Python expressions that comply with solver syntax based on user's query.
}\\
\end{flushleft}
\textbf{Prompt template:}
\begin{flushleft}
\setlength{\rightskip}{0pt plus 1fil}
\raggedright
\texttt{
You are a researcher and have a good knowledge of mathematics and operations research. Your task is to generate a Python expression based on the given notations, ultimately ensuring that the expression meets the syntax requirements of the solver (such as Gurobi) and the output format meets specified standards. \\
For sets, parameters, and variables, the definitions are as follows:\\
\textit{\{Notations\}}\\
Based on the above notations, you are required to provide a Python expression to depict the target:\\
\textit{\{Target\}}\\
The Python expression must be reasonable, containing a part of the mentioned notations without adding other elements, and the formula must be convex. Now I will provide you with some ground-truth formulations for reference.\\
\textit{\{Ground truth\}}\\
The output must contain the following points:\\
- First, provide explanation and illustrate why the formulation is set up in this way;\\
- Second, provide the corresponding Python code for the formula, and only include the given notations;\\
- Third, rewrite Python code to conform to solver syntax.\\
It is necessary to ensure that the output is in the specified format. Take a deep breath and solve the problem step by step. Now you can start working.
}
\end{flushleft}

\subsection{Problem Tailor}
The purpose of the Problem tailor is to identify and exclude parking spots from database tables that are not related to the user's concerns. Furthermore, the Problem tailor provides reasons for selected spots and include contextual information about them. The corresponding prompts are listed as follows: \\
\textbf{Function description:}
\begin{flushleft}
\setlength{\rightskip}{0pt plus 1fil}
\raggedright
\texttt{
Selecting low correlation stations parameterization based on SQL.
}
\end{flushleft}
\textbf{Prompt template:}
\begin{flushleft}
\setlength{\rightskip}{0pt plus 1fil}
\raggedright
\texttt{
You are a database administrator and an expert in operations research. You need to perform cross-analysis on three existing tables in the database and ultimately provide me with some station names and why you selected these station. \\
Explanation of the data tables as follows:\\
\textit{\{Database description\}}\\
Based on the database information provided, you need to tell me which stations have little impact on stations within the target:\\
\textit{\{Target\}}\\
That means it's unnecessary to make decisions, and then all decision variables involved can be parameterized from decisions to parameters. The current Python expression of the target is: \\
\textit{\{Formulation\}}\\
For sets, parameters, and variables, the definitions are as follows:\\
\textit{\{Notations\}}\\
You should pay attention to the following points:\\
- First, you must select at least \textit{\{Min\}} stations and a maximum of \textit{\{Max\}} stations. That means no more than the limit and no less than the requirement. Anyway, you need to give me an answer and you must not choose the stations involved in user targets. \\
- Second, provide me with all the stations to be parameterized for subsequent conversion into parameters. \\
- Third, Explain the reason for your choice.\\
It is necessary to ensure that the output is in the specified format. Take a deep breath and solve the problem step by step. Now you can start working.
}
\end{flushleft}
\subsection{Code Safeguard}
The purpose of the Code safeguard is to go over the erroneous Python code or invalid Gurobi syntax generated by previous steps. After LLM provides an incorrect formulation, Code safeguard will debug each error until the code is able to execute without any issues.  The corresponding prompts are listed as follows:\\
\textbf{Function description:}
\begin{flushleft}
\setlength{\rightskip}{0pt plus 1fil}
\raggedright
\texttt{
Fix errors in the code without changing semantics to ensure it runs properly
}
\end{flushleft}
\textbf{Prompt template:}
\begin{flushleft}
\setlength{\rightskip}{0pt plus 1fil}
\raggedright
\texttt{
You are a programmer specializing in operations research, meticulous and proficient in Gurobi and Python syntax. You are now required to correct the erroneous code. The current Gurobi code in Python is: \\
\textit{\{Gurobi code\}}\\
But this code has some errors, and the detailed error information are as follows:\\
\textit{\{Error traceback\}}\\
You need to modify this code to address the errors specifically, ensuring it can run properly without altering its overall structure. The output must comply with the format requirements.\\
Now I will give you some ground truth formulation. They can be run correctly. If the ground truth I provide to you is valuable, you can refer to it.\\
\textit{\{Ground truth\}}\\
For sets, parameters, and variables, the definitions are as follows: \\
\textit{\{Notations\}}\\
It is necessary to ensure that the output is in the specified format. Take a deep breath and solve the problem step by step. Now you can start working.
}
\end{flushleft}
\subsection{Sample Safeguard}
The purpose of the Sample safeguard is to assess the rationality and compliance of the selected parking spots. After the Problem tailor outputs certain spots, Sample safeguard will reevaluate them and make adjustments if they are deemed inappropriate. The corresponding prompts are listed as follows:\\
\textbf{Function description:}
\begin{flushleft}
\setlength{\rightskip}{0pt plus 1fil}
\raggedright
\texttt{
Check the robustness of the selection to ensure that the results are compliant and reliable.
}
\end{flushleft}
\textbf{Prompt template:}
\begin{flushleft}
\setlength{\rightskip}{0pt plus 1fil}
\raggedright
\texttt{
The selection content given by the previous work is as follows:\\
\textit{\{Selection\}}\\
The selection content is based on a database search to select station for parameterization.\\
You should pay attention to the following points:\\
- Ensure that your output meets the output format specifications. \\
- Verify the reliability of the current selection. The station must be within the given range and the reason for choosing should be sufficient.\\
Verify the above two points, and if they do not meet the requirements, make corrections to your results until they fulfill the request. Now, I'll give you some references to support you in completing tasks.\\
\textit{\{Notations\}}\\
\textit{\{Database description\}}\\
Take a deep breath and solve the problem step by step. Now you can start working.
}
\end{flushleft}

\section{Performance of the Scope-down Policy}\label{app:scope}

We conduct additional experiments to examine how the City-LEO scopes down the EBS parking spots that are relevant to the user's query. 
More precisely, we set the "maximum allowable number of parameterized parking spots" option under the same settings in the aforementioned experiments.
Figure \ref{Fig.max-sites} illustrates the number of parameterized parking spots as the maximum allowed number of parameterized parking spots increases. 
The figure plots both the mean and median numbers of scoped-down parking spots, depicted by a red line for the mean and a black line for the median across different settings of parameterized parking spots (10, 25, 40, 55). The box plots provide an overview of the distribution of scoped-down parking spots at each parameter setting.
\renewcommand{\thefigure}{8}
\begin{figure}[H]
\centering 
\includegraphics[width=0.6\textwidth]{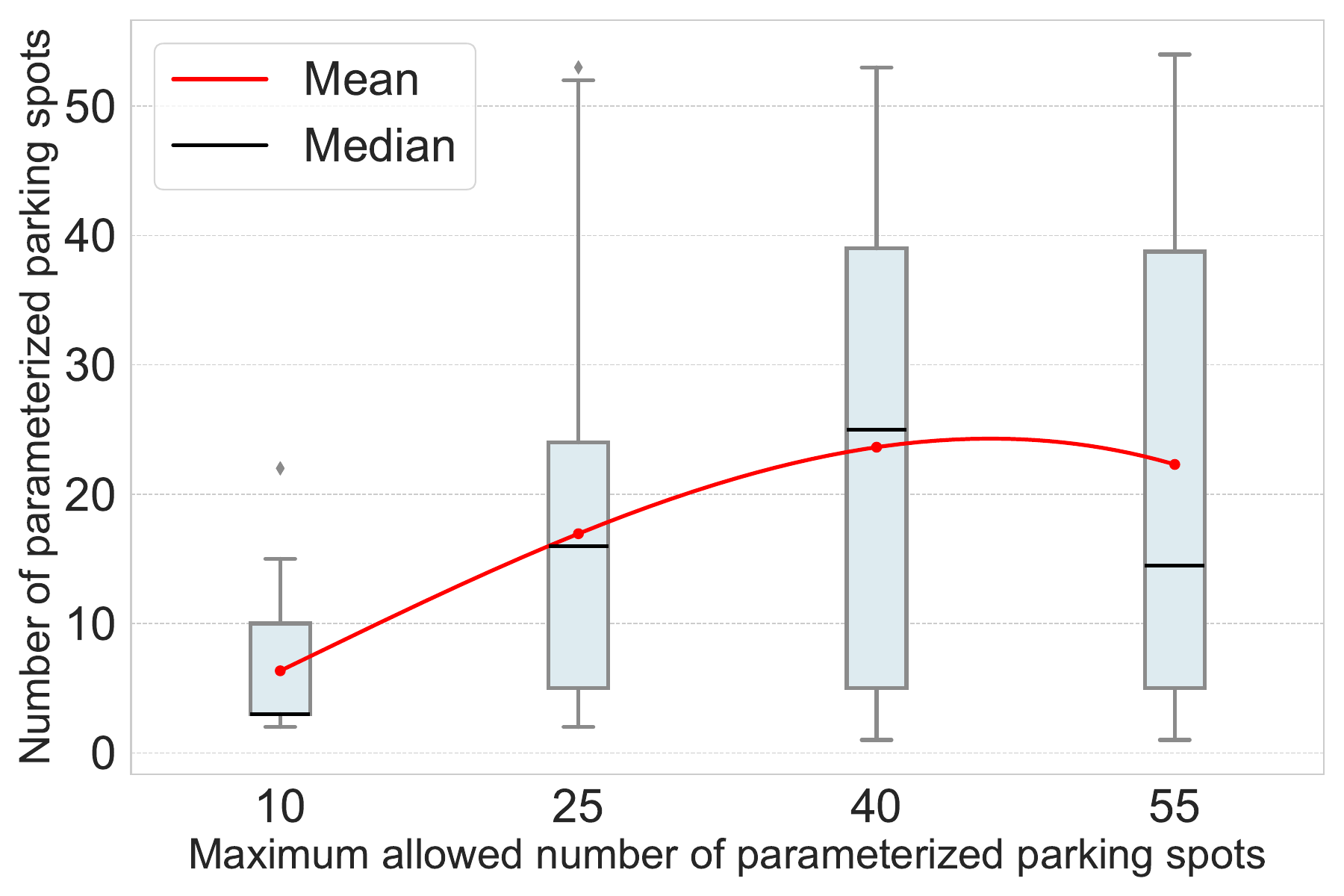} 
\caption{Average Scoped-down Parking Spots Evaluation}{}
 \label{Fig.max-sites} 
\end{figure}
From the Figure \ref{Fig.max-sites}, it can be observed that as the maximum allowed number of parameterized parking spots increases, there is a general trend of increasing mean and median numbers of parameterized parking spots. Specifically, when the maximum allowed parking spots are set at 10, the median number of scoped-down parking spots is relatively low, around 10. This increases to approximately 20, 30, and 40 as the allowed parameterized parking spots are increased to 25, 40, and 55, respectively.

The red line indicating the mean number of scoped-down parking spots shows a steady rise, reinforcing the trend observed in the medians. The spread of the data, as indicated by the whiskers of the box plots, suggests increasing variability in the number of scoped-down parking spots as more spots are allowed to be parameterized, particularly noticeable at the higher settings of 40 and 55. 

\end{document}